\documentclass[12pt]{iopart}
\usepackage{amssymb}
\usepackage{iopams}
\usepackage{graphicx}

\eqnobysec

\begin{document}

\title[Convexification For CIP of MFG]{Convexification for a Coefficient
Inverse Problem of Mean Field Games}
\author{Michael V. Klibanov$^{1}$, Jingzhi Li$^{2}$ and Zhipeng Yang$%
^{3}$}

\begin{abstract}
The globally convergent convexification numerical method is constructed for
a Coefficient Inverse Problem for the Mean Field Games System. A coefficient
characterizing the global interaction term is recovered from the single
measurement data. In particular, a new Carleman estimate for the Volterra
integral operator is proven, and it stronger than the previously known one.
Numerical results demonstrate accurate reconstructions from noisy data.
\end{abstract}

\address{$^1$ Department of Mathematics and Statistics, University of North
Carolina at Charlotte, Charlotte, NC, 28223, USA} 
\address{$^2$ Department
of Mathematics \& National Center for Applied Mathematics Shenzhen \&
SUSTech International Center for Mathematics, Southern University of Science
and Technology, Shenzhen 518055, P.~R.~China} 
\address{$^3$ Department of
Mathematics, Southern University of Science and Technology, Shenzhen 518055,
P.~R.~China} 
\eads{\mailto{mklibanv@uncc.edu},
\mailto{li.jz@sustech.edu.cn}, \mailto{yangzp@sustech.edu.cn}}

\vspace{10pt} \begin{indented}
\item[] 
\end{indented}


\submitto{\IP}

%
%

\section{Introduction}

\label{sec:1}

The mean field games (MFG) theory examines the collective behavior of an
infinite number of rational agents. This theory was initially introduced in
the seminal works of Lasry and Lions \cite{LL1,LL2} and Huang, Caines, and
Malham\'{e} \cite{Huang2,Huang1}. A commonly recognized serious\
mathematical advantage of the MFG theory is that it is based on a universal
system of coupled nonlinear parabolic Partial Differential Equations (PDEs)
known as the Mean Field Games System (MFGS) \cite{A}. At this point of time
the MFG theory is the single mathematical model of a broad spectrum of the
societal phenomena that uses a universal system of PDEs \cite{Burger}.

In addition to the nonlinearity, there are two other quite substantial
challenges in working with the MFGS. The first challenge is that those two
equations have two opposite directions of time. Therefore, the conventional
theory of parabolic PDEs is inapplicable to the MFGS. The second challenge
is due to the presence of an integral operator, the so-called
\textquotedblleft global interaction term", in one of equations of the MFGS,
see below in this section. This presence is very unusual in the theory of
parabolic equations. On one hand, there is no applied meaning of the MFGS
without this term. On the other hand, the presence of this integral term
prevents a straightforward application of previously developed theory of
CIPs for one parabolic PDE to CIPs for the MFGS.

With the increasing significance of social sciences in the modern society,
MFG-based mathematical modeling of social phenomena has a potential of a
substantial societal impact \cite{Burger}. Indeed, this theory finds a
rapidly growing number of applications in such areas as, e.g. finance,
combating corruption, cybersecurity, interactions of electrical vehicles,
election dynamics, robotic control, etc. For a far non-exhaustive list of
references describing these applications, we refer to, e.g. \cite%
{A,Chow,Co,Huang2,Huang1,Kol,LL1,LL2,Trusov}.

\textbf{Definition 1.1}. \emph{We call a problem for the MFGS the
\textquotedblleft forward problem" if it consists in finding a solution of
this system under the assumption that its coefficients are known. And we
call a problem for the MFGS the \textquotedblleft Coefficient Inverse
Problem (CIP)" if it is required in it to find a coefficient(s) of this
system using data occurring in some\ measurement events. }

\textbf{Definition 1.2}. \emph{We call a numerical method for a CIP
\textquotedblleft globally convergent" is there is a rigorous guarantee that
it converges to the true solution of this problem without an advanced
knowledge of any point in a sufficiently small neighborhood of this solution.%
}

Given a broad range of applications of the MFG theory, it becomes important
to address several mathematical questions for both forward problems and CIPs
for the MFGS. We introduce here the \underline{first globally convergent}
numerical method for a CIP for the MFGS. This is a version of the so-called
\textquotedblleft convexification method". The convexification was first
originated in the theoretical work of Klibanov \cite{Klib97}. Next,
publication \cite{Bak} has removed some obstacles for computations, which
were present in \cite{Klib97}. Since then a number of works were published,
which synthesize analytical and numerical parts of various versions of the
convexification for a number of CIPs, see, e.g. \cite{SAR,KL,KLpar,LeLoc}
and references cited therein. In particular, the most recent publication 
\cite{MFG7} of this research group is about the convexification method for a
forward problem for the MFGS.

Conventional \ numerical methods for CIPs are based on the least squares
minimization, see, e.g. \cite{B1,B2,B3,B4,Chavent,Giorgi,Gonch1,Gonch2,Riz}.
While quite powerful, this technique has a drawback: it suffers from the
phenomenon of multiple local minima and ravines of cost functionals, see,
e.g. \cite{Scales} for a good numerical example. This phenomenon, in turn
means that convergence of corresponding algorithms to true solutions can be
rigorously guaranteed only if the starting point of such an algorithm is
located in a sufficiently small neighborhood of that solution. The latter
means \underline{local convergence}. In fact, the convexification method
avoids the phenomenon of local minima and ravines.

This paper consists of two parts:

\begin{enumerate}
\item Construction and global convergence analysis of the convexification
method for our CIP. As a by-product, we establish uniqueness result for our
CIP.

\item Numerical studies. As the second by-product, we present here a new
procedure of data generation for the MFGS.
\end{enumerate}

Another new element of this publication, which is interesting in its own
right, is a new Carleman estimate for a Volterra-like integral operator, see
Theorem 4.1. This estimate is stronger than the previously known one of,
e.g. \cite[Lemma 1.10.3]{KT} and \cite[Lemma 3.1.1]{KL}.

Our CIP is about the recovery of the coefficient $k\left( x\right) $ in 
\begin{equation}
k\left( x\right) \int\limits_{\Omega }K\left( x,y\right) p\left( y,t\right)
dy,  \label{1.1}
\end{equation}%
see section 2 for notations. The integral in (\ref{1.1}) is called the
\textquotedblleft global interaction term". Such a term is not a part of all
previous works on CIPs. The applied meaning of this term is explained in 
\cite[section 5]{MFG1}. More, precisely, the kernel $K(x,y)$ is the
influence on the agent, who occupies the state $x,$ by the agent occupying
the state $y$. The integral (\ref{1.1}) is the average influence on the
agent occupying the state $x,$ by the rest of agents.

Let a small number $\delta >0$ be the level of noise in the input data for
our CIP. To specify Definition 1.2 for our case, the term \textquotedblleft
global convergence" refers to the convergence analysis resulting in a
theorem, which guarantees that if $\delta \rightarrow 0,$ then the iterative
solutions generated by our method converge to the true solution of our CIP
(if it exists) starting from any point within a predefined convex bounded
set in a Hilbert space with a fixed but arbitrary diameter $d>0$. In simple
terms, we do not need a good first guess for the solution. An explicit
estimate of the convergence rate is also given here. Results of our
computational experiments demonstrate a good accuracy of computed solutions
in the presence of a random noise in the input data.

The authors are aware about only two previous publications on numerical
studies of CIPs for the MFGS \cite{Chow,Ding}. Numerical methods of these
references are significantly different from ours, and their global
convergence property is not proven.

Just as in all other previous publications of this research group about both
forward and inverse problems for the MFGS \cite{MFG1}-\cite{MFG7}, we work
here with the input data resulting from a single measurement event. We refer
to \cite{MFG5,MFG6} for previous analytical works on CIPs for the MFGS,
where questions of uniqueness and stability of those problems are addressed
for the single measurement case. In publications \cite{Liu2,Liu1},
uniqueness of some inverse problems for the MFGS was proven for the case of
infinitely many measurements.

The work of Klibanov and Averboukh \cite{MFG1} is the first one, in which
Carleman estimates were introduced in the MFG theory. We quite essentially
use Carleman estimates here, so as in all our works on the MFG theory \cite%
{MFG1}-\cite{MFG7}. Both the convexification method and \cite{MFG5,MFG6} use
the idea of the paper of Bukhgeim and Klibanov \cite{BukhKlib}, where the
method of Carleman estimates was introduced in the field of Inverse Problems
for the first time, see, e.g. \cite{Ksurvey,KL,KT} and references cited
therein for some follow up publications.

\textbf{Remark 1.1}. \emph{Minimal smoothness requirements traditionally
hold little significance in the theory of Ill-Posed and Inverse Problems, as
seen in publications like, e.g. \cite{KL,Nov}, \cite[Theorem 4.1]{Rom}.
Consequently, such requirements are not a primary concern below.}

This paper is arranged as follows. In section 2 we present the MFGS and
formulate our CIP. In section 3 we present the version of the
convexification method for this CIP. In section 4 we formulate theorems of
our global convergence analysis and prove one of them. In sections 5-7 we
prove three more theorems formulated in section 4. Section 8 is dedicated to
numerical studies. Summary of results is given in section 9. All functions
considered below are real valued ones.

\section{Problem Statement}

\label{sec:2}

Let $x=\left( x_{1},...,x_{n}\right) $ denotes points in $\mathbb{R}^{n}.$
Denote $\overline{x}=\left( x_{2},...,x_{n}\right) .$ Let $%
a,b,M,A_{i}>0,i=2,...,n$ and $T>0$ be some numbers and $a<b$. Let the number 
$\gamma \in \left( 0,1\right) .$We set the domain $\Omega \subset \mathbb{R}%
^{n}$ as a rectangular prism, 
\begin{equation}
\left. 
\begin{array}{c}
\Omega =\left\{ x:a<x_{1}<b,-A_{i}<x_{i}<A_{i},i=2,...,n\right\} ,\mbox{ }
\\ 
\Omega _{1}=\left\{ -A_{i}<x_{i}<A_{i},i=2,...,n\right\} ,\mbox{ }\Gamma
=\partial \Omega \cap \left\{ x_{1}=b\right\} , \\ 
Q_{T}=\Omega \times \left( 0,T\right) ,\mbox{ }S_{T}=\partial \Omega \times
\left( 0,T\right) ,\mbox{ }\Gamma _{T}=\Gamma \times \left( 0,T\right) , \\ 
Q_{\gamma T}=\Omega \times \left( \left( 1-\gamma \right) T/2,\left(
1+\gamma \right) T/2\right) \subset Q_{T}.%
\end{array}%
\right.  \label{2.1}
\end{equation}%
We now specify the kernel $K\left( x,y\right) ,$ $x,y\in \Omega $ of the
integral operator (\ref{1.1}) of the global interaction term of the MFGS. It
was noticed in \cite[section 4.2]{LiuOsher} that a good choice for $K\left(
x,y\right) $ would be the product of Gaussians. Since Gaussian approximates
the $\delta -$function in the sense of distributions, then we set%
\begin{equation}
\left. 
\begin{array}{c}
K_{1}\left( x,y\right) =k\left( x\right) \delta \left( x_{1}-y_{1}\right) 
\overline{K}_{1}\left( x,\overline{y}\right) , \\ 
k\left( x\right) \in L_{\infty }\left( \Omega \right) ,\mbox{ }\left\Vert
k\right\Vert _{L_{\infty }\left( \Omega \right) }\leq M, \\ 
\overline{K}_{1}\left( x,\overline{y}\right) \in L_{\infty }\left( \Omega
\times \Omega _{1}\right) ,\mbox{ }\left\Vert \overline{K}_{1}\right\Vert
_{L_{\infty }\left( \Omega \times \Omega _{1}\right) }\leq M.%
\end{array}%
\right.  \label{2.2}
\end{equation}%
We also consider the second form of the kernel $K\left( x,y\right) $ as the
one in \cite{MFG6}:%
\begin{equation}
K_{2}\left( x,y\right) =H\left( y_{1}-x_{1}\right) \overline{K}_{2}\left(
x,y\right) .  \label{2.3}
\end{equation}%
where $H$ is the Heaviside function, 
\[
H\left( z\right) =\left\{ 
\begin{array}{c}
1,z>0, \\ 
0,z<0.%
\end{array}%
\right. 
\]%
Below $K\left( x,y\right) $ is any of functions (\ref{2.2}), (\ref{2.3}).
Let the coefficient 
\begin{equation}
k\left( x\right) \in C\left( \overline{\Omega }\right) .  \label{2.30}
\end{equation}%
Thus, below%
\begin{equation}
\left. 
\begin{array}{c}
k\left( x\right) \int\limits_{\Omega }K\left( x,y\right) f\left( y,t\right)
dy=k\left( x\right) \int\limits_{\Omega _{1}}\overline{K}_{1}\left( x,%
\overline{y}\right) f\left( x_{1},\overline{y},t\right) dy\mbox{,} \\ 
\mbox{in the case }(\ref{2.2}), \\ 
\mbox{and }k\left( x\right) \int\limits_{\Omega }K\left( x,y\right) f\left(
y,t\right) dy=k\left( x\right) \int\limits_{x_{1}}^{b}\left(
\int\limits_{\Omega _{1}}\overline{K}_{2}\left( x,\overline{y}\right)
f\left( y,t\right) d\overline{y}\mbox{ }\right) dy_{1}\mbox{ } \\ 
\mbox{in the case }(\ref{2.3}), \\ 
\forall f\in L_{\infty }\left( Q_{T}\right) .%
\end{array}%
\right.  \label{2.4}
\end{equation}

We consider the MFGS of the second order in the following form \cite{A}:

\begin{equation}
\left. 
\begin{array}{c}
v_{t}(x,t)+\Delta v(x,t){-r(x)(\nabla v(x,t))^{2}/2}- \\ 
-k\left( x\right) \int\limits_{\Omega }K\left( x,y\right) p\left( y,t\right)
dy-s\left( x,t\right) p\left( x,t\right) =0,\mbox{ }\left( x,t\right) \in
Q_{T}, \\ 
p_{t}(x,t)-\Delta p(x,t){-\mbox{div}(r(x)p(x,t)\nabla v(x,t))}=0,\mbox{ }%
\left( x,t\right) \in Q_{T},%
\end{array}%
\right.  \label{2.5}
\end{equation}%
where $s\left( x,t\right) ,p\left( x,t\right) $ is the local interaction
term. Here $v(x,t)$ is the value function and $p(x,t)$ is the density of
players. We assume that%
\begin{equation}
\left. 
\begin{array}{c}
r\in C^{1}\left( \overline{\Omega }\right) ;\mbox{ }s,s_{t}\in L_{\infty
}\left( Q_{T}\right) , \\ 
\left\Vert r\right\Vert _{C^{1}\left( \overline{\Omega }\right) },\left\Vert
s\right\Vert _{L_{\infty }\left( Q_{T}\right) },\left\Vert s_{t}\right\Vert
_{L_{\infty }\left( Q_{T}\right) }\leq M.%
\end{array}%
\right.  \label{2.40}
\end{equation}

\textbf{Coefficient Inverse Problem} (CIP). \emph{Assume that the following
functions }$v_{0}\left( x\right) $, $p_{0}\left( x\right) $, $g_{0,i}\left(
x,t\right) $, $g_{1,i}\left( x,t\right) $\emph{\ are known:}%
\begin{equation}
\left. 
\begin{array}{c}
v\left( x,T/2\right) =v_{0}\left( x\right) ,\mbox{ }p\left( x,T/2\right)
=p_{0}\left( x\right) ,\mbox{ }x\in \Omega , \\ 
v\mid _{S_{T}}=g_{0,1}\left( x,t\right) ,\mbox{ }p\mid
_{S_{T}}=g_{0,2}\left( x,t\right) , \\ 
v_{x_{1}}\mid _{\Gamma _{T}}=g_{1,1}\left( x,t\right) ,\mbox{ }p_{x_{1}}\mid
_{\Gamma _{T}}=g_{1,2}\left( x,t\right) .%
\end{array}
\right.  \label{2.400}
\end{equation}%
\emph{Determine the coefficient }$k\left( x\right) $\emph{\ in (\ref{2.5}),
assuming (\ref{2.30}).}

Note that since the Neumann boundary conditions are not given at $%
S_{T}\diagdown \Gamma _{T}$ then the boundary data in (\ref{2.5}) are
incomplete. In (\ref{2.5})\emph{\ }$T/2$ can be replaced with any fixed
number $t_{0}\in \left( 0,T\right) $ without any significant changes in our
derivations below.

\textbf{Remark 2.1}. \emph{The data (\ref{2.5}) are generated by a single
measurement event. As to a real game, it was pointed out in \cite[section
2.2 ]{MFG6} that these data can be obtained via polling game participants.}

\section{Convexification}

\label{sec:3}

Below we work either with (\ref{2.2}) or with (\ref{2.3}) and keep (\ref{2.4}%
).\emph{\ }Let $c>0$ be a number. Assume that 
\begin{equation}
\mbox{in the case }\mbox{(\ref{2.2})}\mbox{ }\left\vert \int\limits_{\Omega
_{1}}\overline{K}_{1}\left( x,\overline{y}\right) p_{0}\left( x_{1},%
\overline{y}\right) d\overline{y}\right\vert \geq c,\mbox{ }x\in \overline{%
\Omega }.  \label{3.1}
\end{equation}%
Let $\eta \in \left( 0,b-a\right) $ be a number. Denote%
\[
\chi \left( x_{1}\right) =\left\{ 
\begin{array}{c}
1,x_{1}\in \left[ 0,b-\eta \right) , \\ 
0,x_{1}\in \left[ b-\eta ,b\right] .%
\end{array}%
\right\} 
\]%
In addition to (\ref{3.1}), we also assume that  
\begin{equation}
\left. 
\begin{array}{c}
\mbox{In the case (\ref{2.3})} \\ 
k\left( x\right) =0,\mbox{ }x_{1}\in \left( b-\eta ,b\right) , \\ 
\chi \left( x_{1}\right) \int\limits_{x_{1}}^{a}\left( \int\limits_{\Omega
_{1}}\overline{K}_{2}\left( x,\overline{y}\right) p_{0}\left( y\right) d%
\overline{y}\mbox{ }\right) dy_{1}\geq c,\mbox{ }x\in \overline{\Omega }.%
\end{array}%
\right.   \label{3.01}
\end{equation}%
Denote%
\begin{equation}
u\left( x,t\right) =v_{t}\left( x,t\right) ,\mbox{ }m\left( x,t\right)
=p_{t}\left( x,t\right) .  \label{3.2}
\end{equation}%
Using (\ref{2.400}) and (\ref{3.2}), we obtain%
\begin{equation}
\left. 
\begin{array}{c}
\hspace{-1.8cm}v(x,t)=\int\limits_{T/2}^{t}u(x,\tau )d\tau +v_{0}\left(
x\right) ,\mbox{ }\left( x,t\right) \in Q_{T}, \\ 
p(x,t)=\int\limits_{T/2}^{t}m(x,\tau )d\tau +p_{0}\left( x\right) ,\mbox{ }%
\left( x,t\right) \in Q_{T}.%
\end{array}%
\right.   \label{3.3}
\end{equation}%
We have:%
\begin{equation}
u\left( x,T/2\right) =u\left( x,t\right) -\int\limits_{T/2}^{t}u_{t}\left(
x,\tau \right) d\tau .  \label{3.4}
\end{equation}%
Setting in the first line of (\ref{2.5}) $t=T/2$ and using (\ref{3.1})-(\ref%
{3.4}), we obtain%
\begin{equation}
\left. 
\begin{array}{c}
k\left( x\right) =u\left( x,T/2\right) f\left( x\right) +F\left( x\right) ,
\\ 
\mbox{or, equivalently, }k\left( x\right) =\left( u\left( x,t\right)
-\int\limits_{T/2}^{t}u_{t}\left( x,\tau \right) d\tau \right) f\left(
x\right) +F\left( x\right) ,\mbox{ } \\ 
F\left( x\right) =\left[ \left( \Delta v_{0}-r\left( \nabla v_{0}\right)
^{2}/2-s\left( x,T/2\right) p_{0}\right) \left( x\right) \right] f\left(
x\right) ,%
\end{array}%
\right.   \label{3.5}
\end{equation}%
where the function $f\left( x\right) $ is defined as: 
\begin{equation}
\left. f\left( x\right) =\left\{ 
\begin{array}{c}
\left[ \int\limits_{\Omega _{1}}\overline{K}_{1}\left( x,\overline{y}\right)
p_{0}\left( x_{1},\overline{y}\right) d\overline{y}\right] ^{-1}, \\ 
\mbox{in the case }\emph{(\ref{2.2}),} \\ 
\chi \left( x_{1}\right) \left[ \int\limits_{x_{1}}^{a}\left(
\int\limits_{\Omega _{1}}\overline{K}_{2}\left( x,\overline{y}\right)
p_{0}\left( y\right) d\overline{y}\mbox{ }\right) dy_{1}\right] ^{-1}, \\ 
\mbox{in the case }\emph{(\ref{2.3}).}%
\end{array}%
\right\} \right.   \label{3.6}
\end{equation}

Differentiate equations (\ref{2.5}) with respect to $t.$ Substituting (\ref%
{3.2})-(\ref{3.6}) in resulting equations and using (\ref{2.4}), we obtain
two nonlinear integral differential equations with the lateral Cauchy data
for the vector function $\left( u,m\right) \left( x,t\right) .$ The first
equation is:%
\begin{equation}
\hspace{-1 cm}
\left. 
\begin{array}{c}
L_{1}\left( u,m,v_{0},p_{0}\right) =u_{t}+\Delta u{-r\nabla u}\left(
\int\limits_{T/2}^{t}\nabla u(x,\tau )d\tau +\nabla v_{0}\left( x\right)
\right) - \\ 
-f\left( x\right) \left[ \left( u\left( x,t\right)
-\int\limits_{T/2}^{t}u_{t}\left( x,\tau \right) d\tau \right) +F\left(
x\right) \right] \int\limits_{\Omega }K\left( x,y\right) m\left( y,t\right)
dy- \\ 
-sm-s_{t}\left( \int\limits_{T/2}^{t}m(x,\tau )d\tau +p_{0}\left( x\right)
\right) ,\mbox{ }\left( x,t\right) \in Q_{T}.%
\end{array}
\right.  \label{3.61}
\end{equation}%
The second equation is: 
\begin{equation}
\left. 
\begin{array}{c}
L_{2}\left( u,m,v_{0},p_{0}\right) =m_{t}-\Delta m(x,t){-} \\ 
-{\mbox{div}}\left[ {r(x)m\left( \int\limits_{T/2}^{t}\nabla u(x,\tau )d\tau
+\nabla v_{0}\left( x\right) \right) }\right] - \\ 
-\mbox{div}\left[ {r(x)\nabla u}\left( \int\limits_{T/2}^{t}m(x,\tau )d\tau
+p_{0}\left( x\right) \right) \right] =0,\mbox{ }\left( x,t\right) \in Q_{T}.%
\end{array}
\right.  \label{3.62}
\end{equation}%
The lateral Cauchy data for functions $u\left( x,t\right) $ and $m\left(
x,t\right) $ are:%
\begin{equation}
\left. 
\begin{array}{c}
u\mid _{S_{T}}=\partial _{t}g_{0,1}\left( x,t\right) ,\mbox{ }m\mid
_{S_{T}}=\partial _{t}g_{0,2}\left( x,t\right) , \\ 
u_{x_{1}}\mid _{\Gamma _{T}}=\partial _{t}g_{1,1}\left( x,t\right) ,\mbox{ }
m_{x_{1}}\mid _{\Gamma _{T}}=\partial _{t}g_{1,2}\left( x,t\right) .%
\end{array}
\right.  \label{3.63}
\end{equation}

Suppose that problem (\ref{3.61})-(\ref{3.63}) is solved. Then the unknown
coefficient $k\left( x\right) $ can be recovered via the first line of (\ref%
{3.5}). Therefore, we focus below on the numerical solution of problem (\ref%
{3.61})-(\ref{3.63}). In our derivations below we need \ functions $%
u_{x_{i}},u_{t},m_{x_{i}},m_{t},\Delta u,\Delta m\in C\left( \overline{Q}%
_{T}\right) ,$ see Remark 1.1. Therefore, we assume that functions $u,m\in
H^{k_{n}}\left( Q_{T}\right) ,$ where $k_{n}=\left[ \left( n+1\right) /2%
\right] +3,$ where $\left[ \left( n+1\right) /2\right] $ is the largest
integer not exceeding $\left( n+1\right) /2.$ By embedding theorem there
exists a constant $C_{0}=C_{0}\left( Q_{T}\right) >0$ depending only on the
domain $Q_{T}$ such that%
\begin{equation}
\hspace{-1cm}H^{k_{n}}\left( Q_{T}\right) \subset C^{2}\left( \overline{Q}%
_{T}\right) ,\mbox{ }\left\Vert y\right\Vert _{C^{2}\left( \overline{Q}%
_{T}\right) }\leq C_{0}\left\Vert y\right\Vert _{H^{k_{n}}\left(
Q_{T}\right) },\mbox{ }\forall y\in H^{k_{n}}\left( Q_{T}\right) .
\label{3.7}
\end{equation}

Introduce four Hilbert spaces: 
\begin{equation}
\left. 
\begin{array}{c}
H^{2,1}\left( Q_{T}\right) =\left\{ u:\left\Vert u\right\Vert
_{H^{2,1}\left( Q_{T}\right) }^{2}=\sum\limits_{\left\vert s\right\vert
+2k\leq 2}\left\Vert D_{x}^{s}D_{t}^{k}u\right\Vert _{H^{2,1}\left(
Q_{T}\right) }^{2}<\infty \right\} , \\ 
H_{0}^{2,1}\left( Q_{T}\right) =\left\{ u\in H^{2,1}\left( Q_{T}\right)
:u\mid _{S_{T}}=u_{x_{1}}\mid _{\Gamma _{T}}=0\right\} , \\ 
H=\left\{ \left( u,m\right) :\left\Vert \left( u,m\right) \right\Vert
_{H}^{2}=\left\Vert u\right\Vert _{H^{k_{n}}\left( Q_{T}\right)
}^{2}+\left\Vert m\right\Vert _{H^{k_{n}}\left( Q_{T}\right) }^{2}<\infty
\right\} , \\ 
H_{0}=\left\{ 
\begin{array}{c}
\left( u,m\right) \in H: \\ 
u\mid _{S_{T}}=u_{x_{1}}\mid _{\Gamma _{T}}=m\mid _{S_{T}}=m_{x_{1}}\mid
_{\Gamma _{T}}=0.%
\end{array}%
\right\}%
\end{array}%
\right.  \label{4.01}
\end{equation}

\textbf{Remark 3.1}. \emph{Below }$\left[ ,\right] $\emph{\ is the scalar
product in the space }$H.$

Let $R>0$ be an arbitrary number. Consider two sets $B\left( R\right) $ and $%
B_{0}\left( R\right) $ defined as:%
\begin{equation}
B\left( R\right) =\left\{ 
\begin{array}{c}
\left( u,m\right) \in H: \\ 
u\mid _{S_{T}}=\partial _{t}g_{0,1}\left( x,t\right) ,\mbox{ }m\mid
_{S_{T}}=\partial _{t}g_{0,2}\left( x,t\right) , \\ 
u_{x_{1}}\mid _{\Gamma _{T}}=\partial _{t}g_{1,1}\left( x,t\right) ,\mbox{ }%
m_{x_{1}}\mid _{\Gamma _{T}}=\partial _{t}g_{1,2}\left( x,t\right) , \\ 
\left\Vert \left( u,m\right) \right\Vert _{H}<R,%
\end{array}%
\right\}  \label{3.8}
\end{equation}%
\begin{equation}
B_{0}\left( R\right) =\left\{ \left( u,m\right) \in H_{0}:\left\Vert \left(
u,m\right) \right\Vert _{H}<R\right\} .  \label{3.80}
\end{equation}%
Let $\alpha $ be a rational number represented as%
\begin{equation}
\alpha =\frac{n_{1}}{n_{2}}\in \left( 0,\frac{1}{3}\right) ,  \label{3.9}
\end{equation}%
where $n_{1}>0$ and $n_{2}>0$ are two odd integers. Let $\lambda \geq 1$ be
a large parameter which will be chosen later. We now consider a new Carleman
Weight Function (CWF) $\varphi _{\lambda }\left( x,t\right) ,$%
\begin{equation}
\varphi _{\lambda }\left( x,t\right) =\exp \left[ 2\lambda \left(
x_{1}^{2}-\left( t-T/2\right) ^{1+\alpha }\right) \right] .  \label{3.10}
\end{equation}%
Note that since $1+\alpha =\left( n_{1}+n_{2}\right) /n_{2}$ and since the
number $n_{1}+n_{2}$ is even, then $\left( t-T/2\right) ^{1+\alpha }$ is
defined for both $t>T/2$ and $t<T/2,$ and the function $\varphi _{\lambda
}\left( x,t+T/2\right) $ is even with respect to $t$. The novelty of the CWF
(\ref{3.10}) is due to the fact that the case $1+\alpha =2$ in the
conventional one, see, e.g. \cite[formula (3.12)]{KLpar} and \cite[formula
(9.20)]{KL}. Below we need the form (\ref{3.10}) with $1+\alpha \in \left(
1,4/3\right) $ in order to obtain proper estimates for the Volterra
integrals in (\ref{3.61}), (\ref{3.62}) containing $u_{t}$ and $\Delta u.$
Clearly 
\begin{equation}
\max_{\overline{Q}_{T}}\varphi _{\lambda }\left( x,t\right) =e^{2\lambda
b^{2}}.  \label{3.100}
\end{equation}

Consider four functionals mapping the set $\overline{B\left( R\right) }$ in $%
\mathbb{R},$%
\begin{equation}
\hspace{-2cm}\left. 
\begin{array}{c}
J_{1,\lambda }\left( u,m\right) =\int\limits_{Q_{T}}\left( L_{1}\left(
u,m\right) \right) ^{2}\varphi _{\lambda }dxdt, \\ 
J_{2,\lambda }\left( u,m\right) =\int\limits_{Q_{T}}\left( L_{2}\left(
u,m\right) \right) ^{2}\varphi _{\lambda }dxdt, \\ 
J_{3,\beta }\left( u,m\right) =\beta \left\Vert \left( u,m\right)
\right\Vert _{H}^{2}, \\ 
J_{\lambda ,\beta }\left( u,m\right) =e^{-2\lambda b^{2}}\lambda
^{3/2}J_{1,\lambda }\left( u,m\right) +e^{-2\lambda b^{2}}J_{2,\lambda
}\left( u,m\right) +J_{3,\beta }\left( u,m\right) ,%
\end{array}%
\right.  \label{3.11}
\end{equation}%
where the operators $L_{1}\left( u,p\right) $ and $L_{2}\left( u,p\right) $
are defined in (\ref{3.61}) and (\ref{3.62}), and $\beta \in \left(
0,1\right) $ is the regularization parameter. The multiplier $e^{-2\lambda
b^{2}}$ is included in $J_{\lambda ,\beta }$ to balance the terms in the
first two lines of (\ref{3.11}) with the term in the third line, see (\ref%
{3.100}). We solve problem (\ref{3.61})-(\ref{3.63}) via solving the
following Minimization Problem:

\textbf{Minimization Problem}. \emph{Minimize the functional }$J_{\lambda
,\beta }\left( u,m\right) $\emph{\ in (\ref{3.11}) on the set }$\overline{
B\left( R\right) }$\emph{\ defined in (\ref{3.8}).}

\section{Convergence Analysis}

\label{sec:4}

\subsection{Carleman estimates}

\label{sec:4.1}

\textbf{Theorem 4.1 }(a new Carleman estimate for a Volterra-like integral). 
\emph{Let }$d>0$\emph{\ be a number and let the number }$\alpha $ \emph{be
as in (\ref{3.9}), where }$n_{1},n_{2}$\emph{\ are two odd numbers. Then the
following Carleman estimate of the Volterra-like integral holds for all
functions }$f\in L_{2}\left( -d,d\right) $ \emph{and for all} $\lambda >0:$ 
\emph{\ }%
\begin{equation}
\int\limits_{-d}^{d}e^{-2\lambda t^{1+\alpha }}\left(
\int\limits_{0}^{t}f\left( \tau \right) d\tau \right) ^{2}dt\leq \frac{1}{%
\lambda ^{3/2}}\cdot \frac{d^{\left( 1-3\alpha \right) /2}}{\sqrt{2}\left(
1+\alpha \right) ^{3/2}}\int\limits_{-d}^{d}f^{2}e^{-2\lambda t^{1+\alpha
}}dt.  \label{4.1}
\end{equation}

\textbf{Remark 4.1.} \emph{As stated in section 3, }$1+\alpha =2$ \emph{in
the conventional case.\ However, the corresponding conventional analog of
estimate (\ref{4.1})\ is\ weaker than (\ref{4.1})\ since\ }$\lambda ^{-3/2}$%
\emph{\ is replaced then with }$\lambda ^{-1}>>\lambda ^{-3/2}$\emph{\ \ for
large values of }$\lambda ,$ \emph{see \cite[Lemma 3.1.1]{KT}, \cite[ Lemma
3.1.1]{KL} for the conventional case. We also refer to the proof of Theorem
4.3 for a more detailed explanation.}

\textbf{Theorem 4.2.}\emph{\ Let }$\varphi _{\lambda }\left( x,t\right) $%
\emph{\ be the Carleman Weight Function defined in (\ref{3.10}). Then there
exist a sufficiently large number }$\lambda _{0}=\lambda _{0}\left( \Omega
,\alpha ,T\right) \geq 1$\emph{\ and a number }$C_{0}=C_{0}\left( \Omega
,\alpha ,T\right) >0,$\emph{\ both numbers depending only on listed
parameters, such that the following two Carleman estimates are valid: }%
\begin{equation}
\hspace{-1cm}\left. 
\begin{array}{c}
\int\limits_{Q_{T}}\left( u_{t}\pm \Delta u\right) ^{2}\varphi _{\lambda
}dxdt\geq \left( C_{0}/\lambda \right) \int\limits_{Q_{T}}\left(
u_{t}^{2}+\sum\limits_{i,j=1}^{n}u_{x_{i}x_{j}}^{2}\right) \varphi _{\lambda
}dxdt+ \\ 
+C_{0}\int\limits_{Q_{T}}\left( \lambda \left( \nabla u\right) ^{2}+\lambda
^{3}u^{2}\right) \varphi _{\lambda }dxdt- \\ 
-C_{0}\left( \left\Vert u\left( x,T\right) \right\Vert _{H^{1}\left( \Omega
\right) }^{2}+\left\Vert u\left( x,0\right) \right\Vert _{H^{1}\left( \Omega
\right) }^{2}\right) \exp \left[ -2\lambda \left( \left( T/2\right)
^{1+\alpha }-b^{2}\right) \right] ,\mbox{ } \\ 
\forall \lambda \geq \lambda _{0},\mbox{ }\forall u\in H_{0}^{2,1}\left(
Q_{T}\right) .%
\end{array}%
\right.  \label{4.2}
\end{equation}

An analog of Theorem 4.2 was proven in \cite{KLpar} and \cite[Theorem 9.4.1]%
{KL} for the case $1+\alpha =2$ and for the operator $\partial _{t}-\Delta .$
The proof for our case is quite similar with some insignificant deviations.
Therefore, we omit the proof of this theorem.

\subsection{Global strong convexity and uniqueness}

\label{sec:4.2}

\textbf{Theorem 4.3} (the central result). \emph{Assume that }%
\begin{equation}
\left\Vert v_{0}\right\Vert _{C^{2}\left( \overline{\Omega }\right) },%
\mbox{ 
	}\left\Vert p_{0}\right\Vert _{C^{1}\left( \overline{\Omega }\right) }<R
\label{4.03}
\end{equation}%
\emph{and conditions (\ref{2.1})-(\ref{2.4}), (\ref{3.1}) and (\ref{3.01})
are satisfied. Then:}

\emph{1. The functional }$J_{\lambda ,\beta }\left( u,m\right) $\emph{\ has
Fr\'{e}chet derivative }$J_{\lambda ,\beta }^{\prime }\left( u,m\right) \in
H_{0}$\emph{\ at every point }$\left( u,m\right) \in \overline{B\left(
R\right) },$\emph{\ and this derivative is Lipschitz continuous on }$%
\overline{B\left( R\right) }$\emph{, i.e. there exists a number }$D=D\left(
\Omega ,T,M,c,\alpha ,R,\lambda \right) >0$\emph{\ depending only on listed
parameters such that for all }$\left( u_{1},m_{1}\right) ,\left(
u_{2},m_{2}\right) \in \overline{B\left( R\right) }$\emph{\ }%
\begin{equation}
\left\Vert J_{\lambda ,\beta }^{\prime }\left( u_{2},m_{2}\right)
-J_{\lambda ,\beta }^{\prime }\left( u_{1},m_{1}\right) \right\Vert _{H}\leq
D\left\Vert \left( u_{1},m_{1}\right) -\left( u_{2},m_{2}\right) \right\Vert
_{H}.  \label{4.3}
\end{equation}

\emph{2. Let }$\lambda _{0}=\lambda _{0}\left( \Omega ,\alpha ,T\right) \geq
1$\emph{\ be the number of Theorem 4.2. There exist a sufficiently large
number }$\lambda _{1}=\lambda _{1}\left( \Omega ,\alpha ,T,M,c,R\right) \geq
\lambda _{0}$\emph{\ and a number }$C_{1}=C_{1}\left( \Omega ,\alpha
,T,M,c,R\right) >0,$\emph{\ both numbers depending only on listed
parameters, such that if }$\lambda \geq \lambda _{1}$\emph{\ and }$\beta \in %
\left[ 2e^{-\lambda \left( T/2\right) ^{1+\alpha }},1\right] ,$\emph{\ then
the functional }$J_{\lambda ,\beta }$\emph{\ is strongly convex on }$%
\overline{B\left( R\right) },$\emph{\ i.e. the following inequality holds:}%
\begin{equation}
\hspace{-1cm}\left. 
\begin{array}{c}
J_{\lambda ,\beta }\left( u_{2},m_{2}\right) -J_{\lambda ,\beta }\left(
u_{1},m_{1}\right) -\left[ J_{\lambda ,\beta }^{\prime }\left(
u_{1},m_{1}\right) ,\left( u_{2},m_{2}\right) -\left( u_{1},m_{1}\right) %
\right] \geq \\ 
\geq C_{1}\exp \left[ -2\lambda \left( \gamma T/2\right) ^{1+\alpha }\right]
\left\Vert \left( u_{2},m_{2}\right) -\left( u_{1},m_{1}\right) \right\Vert
_{H^{2,1}\left( Q_{\gamma T}\right) \times H^{2,1}\left( Q_{\gamma T}\right)
}^{2}+ \\ 
+\left( \beta /2\right) \left\Vert \left( u_{2},m_{2}\right) -\left(
u_{1},m_{1}\right) \right\Vert _{H}^{2},\mbox{ } \\ 
\forall \left( u_{1},m_{1}\right) ,\left( u_{2},m_{2}\right) \in \overline{%
B\left( R\right) },\mbox{ }\forall \geq \lambda _{1},\mbox{ }\forall \gamma
\in \left( 0,1\right) .%
\end{array}%
\right.  \label{4.4}
\end{equation}

\emph{3. For }$\lambda $\emph{\ and }$\beta $\emph{\ as in item 2, the
functional }$J_{\lambda ,\beta }$\emph{\ has unique minimizer }$\left(
u_{\min ,\lambda ,\beta },m_{\min ,\lambda ,\beta }\right) $\emph{\ on the
set }$\overline{B\left( R\right) },$\emph{\ and the following inequality
holds:}%
\begin{equation}
\hspace{-1.5cm}\left[ J_{\lambda ,\gamma }^{\prime }\left( u_{\min ,\lambda
,\beta },m_{\min ,\lambda ,\beta }\right) ,\left( u_{\min ,\lambda ,\beta
}-u,m_{\min ,\lambda ,\beta }-m\right) \right] \leq 0,\mbox{ }\forall \left(
u,m\right) \in \overline{B\left( R\right) }.  \label{4.5}
\end{equation}

\textbf{Remarks 4.2:}

\begin{enumerate}
\item \emph{Below }$C_{1}>0$\emph{\ denotes different numbers depending only
on the above listed parameters.}

\item \emph{Even though this theorem requires that }$\lambda $\emph{\ should
be sufficiently large, our extensive computational experience with the
convexification method tells us that reasonable values of }$\lambda \in %
\left[ 1,5\right] $\emph{\ can always be selected to obtain good accuracy of
numerical results, see, e.g. \cite{Bak,KLpar,KL,SAR,MFG7,LeLoc} and
references cited therein as well as subsection 6.2 below. This is basically
because, like in any asymptotic theory, only numerical studies can indicate
which specific ranges of parameters are reasonable.}
\end{enumerate}

\subsection{The accuracy of the minimizer and uniqueness}

\label{sec:4.3}

Suppose that there exists a 2d vector function $G=\left( G_{1},G_{2}\right)
\in B\left( R\right) .$ By (\ref{3.8}) this means that this vector function
is an extension of boundary conditions of (\ref{3.63}) inside of the domain $%
Q_{T},$%
\begin{equation}
\left. 
\begin{array}{c}
G_{1}\mid _{S_{T}}=\partial _{t}g_{0,1}\left( x,t\right) ,\mbox{ }G_{2}\mid
_{S_{T}}=\partial _{t}g_{0,2}\left( x,t\right) , \\ 
\partial _{x_{1}}G_{1}\mid _{\Gamma _{T}}=\partial _{t}g_{1,1}\left(
x,t\right) ,\mbox{ }\partial _{x_{1}}G_{2}\mid _{\Gamma _{T}}=\partial
_{t}g_{1,2}\left( x,t\right) .%
\end{array}%
\right.  \label{4.6}
\end{equation}%
By one of principles of the theory of Ill-Posed Problems \cite{T}, we assume
that there exists an exact solution $\left( k^{\ast },u^{\ast },m^{\ast
}\right) \in C\left( \overline{\Omega }\right) \times B^{\ast }\left(
R-\delta \right) $ of CIP (\ref{2.30}), (\ref{2.400}) with the exact,
noiseless data 
\begin{equation}
v_{0}^{\ast },p_{0}^{\ast },\partial _{t}g_{0,1}^{\ast },\partial
_{t}g_{0,2}^{\ast },\partial _{t}g_{1,1}^{\ast },\partial _{t}g_{1,2}^{\ast
},  \label{4.7}
\end{equation}%
where%
\begin{equation}
B^{\ast }\left( R-\delta \right) =\left\{ 
\begin{array}{c}
\left( u,m\right) \in H: \\ 
u\mid _{S_{T}}=\partial _{t}g_{0,1}^{\ast },\mbox{ }m\mid _{S_{T}}=\partial
_{t}g_{0,2}^{\ast }, \\ 
u_{x_{1}}\mid _{\Gamma _{T}}=\partial _{t}g_{1,1}^{\ast },\mbox{ }%
m_{x_{1}}\mid _{\Gamma _{T}}=\partial _{t}g_{1,2}^{\ast }, \\ 
\left\Vert \left( u,m\right) \right\Vert _{H}<R-\delta .%
\end{array}%
\right\} ,  \label{4.71}
\end{equation}%
where $\delta >0$ is a sufficiently small number characterizing the level of
noise in the input data. Hence, there exists a vector function $G^{\ast
}=\left( G_{1}^{\ast },G_{2}^{\ast }\right) \in B^{\ast }\left( R\right) $
satisfying direct analogs of boundary conditions (\ref{4.6}), in which the
right hand sides are replaced with functions listed in (\ref{4.7}). Thus,%
\begin{equation}
\left( u^{\ast },m^{\ast }\right) ,G^{\ast }=\left( G_{1}^{\ast
},G_{2}^{\ast }\right) \in B^{\ast }\left( R-\delta \right) .  \label{4.72}
\end{equation}%
Furthermore, by (\ref{3.5}) the exact coefficient $k^{\ast }\left( x\right) $
is 
\begin{equation}
k^{\ast }\left( x\right) =u^{\ast }\left( x,T/2\right) f^{\ast }\left(
x\right) +F^{\ast }\left( x\right) ,  \label{4.70}
\end{equation}%
where functions $f^{\ast }\left( x\right) $ and $F^{\ast }\left( x\right) $
are obtained from functions $f\left( x\right) $ and $F\left( x\right) $ in (%
\ref{3.5}), (\ref{3.6}) via replacing the pair $\left( v_{0},p_{0}\right) $
with the pair $\left( v_{0}^{\ast },p_{0}^{\ast }\right) .$

We assume that 
\begin{equation}
\left. 
\begin{array}{c}
\left\Vert v_{0}^{\ast }\right\Vert _{C^{2}\left( \overline{\Omega }\right)
},\left\Vert p_{0}^{\ast }\right\Vert _{C^{1}\left( \overline{\Omega }
\right) }<R, \\ 
\left\Vert G-G^{\ast }\right\Vert _{H}<\delta , \\ 
\left\Vert v_{0}-v_{0}^{\ast }\right\Vert _{C^{2}\left( \overline{\Omega }
\right) },\left\Vert p_{0}-p_{0}^{\ast }\right\Vert _{C^{1}\left( \overline{
\Omega }\right) }<\delta .%
\end{array}
\right.  \label{4.8}
\end{equation}

Let $\left( u,m\right) \in B\left( R\right) $ be an arbitrary pair of
functions. Denote%
\begin{equation}
\left( \widetilde{u}^{\ast },\widetilde{m}^{\ast }\right) =\left( u^{\ast
},m^{\ast }\right) -G^{\ast },\mbox{ }\left( \widetilde{u},\widetilde{m}%
\right) =\left( u,m\right) -G.  \label{4.9}
\end{equation}%
Note that by (\ref{4.71}), (\ref{4.72}) and (\ref{4.8}) 
\begin{equation}
\left\Vert G\right\Vert _{H}<R.  \label{4.90}
\end{equation}%
By (\ref{4.01})-(\ref{3.80}) and (\ref{4.8})-(\ref{4.90})%
\begin{equation}
\left( \widetilde{u}^{\ast },\widetilde{m}^{\ast }\right) ,\left( \widetilde{%
u},\widetilde{m}\right) \in \overline{B_{0}\left( 2R\right) }.  \label{4.10}
\end{equation}%
Also, let $L_{i}\left( u,m,v_{0},p_{0}\right) ,$ $i=1,2$ be operators
defined in (\ref{3.61}), (\ref{3.62}). Then%
\begin{equation}
L_{i}\left( u^{\ast },m^{\ast },v_{0}^{\ast },p_{0}^{\ast }\right)
=L_{i}\left( \widetilde{u}^{\ast }+G_{1}^{\ast },\widetilde{m}^{\ast
}+G_{2}^{\ast },v_{0}^{\ast },p_{0}^{\ast }\right) =0,\mbox{ }i=1,2.
\label{4.11}
\end{equation}%
Introduce a new functional $I_{\lambda ,\beta },$ 
\begin{equation}
I_{\lambda ,\beta }:\overline{B_{0}\left( 2R\right) }\rightarrow \mathbb{R},%
\mbox{ }I_{\lambda ,\beta }\left( \widetilde{u},\widetilde{m}\right)
=J_{\lambda ,\beta }\left( \widetilde{u}+G_{1},\widetilde{m}+G_{2}\right) .
\label{4.12}
\end{equation}%
It follows from (\ref{3.8}), (\ref{3.80}), (\ref{4.90}) and triangle
inequality that 
\begin{equation}
\left( \widetilde{u}+G_{1},\widetilde{m}+G_{2}\right) \in \overline{B\left(
3R\right) },\mbox{ }\forall \left( \widetilde{u},\widetilde{m}\right) \in 
\overline{B_{0}\left( 2R\right) }.  \label{4.120}
\end{equation}

\textbf{Theorem 4.4} (the accuracy of the minimizer and uniqueness of the
CIP). \emph{\ Suppose that conditions of Theorem 4.3 as well as conditions (%
\ref{4.7})-( \ref{4.12}) are met. In addition, let inequalities (\ref{3.1})
and (\ref{3.01}) be valid when }$p_{0}$\emph{\ is replaced with }$%
p_{0}^{\ast }.$\emph{\ Then:}

\emph{1. The functional }$I_{\lambda ,\beta }$\emph{\ has the Fr\'{e}chet
derivative }$I_{\lambda ,\beta }^{\prime }\left( \widetilde{u},\widetilde{m}%
\right) \in H_{0}$\emph{\ at any point }$\left( \widetilde{u},\widetilde{m}%
\right) \in \overline{B_{0}\left( 2R\right) }$\emph{\ and the analog of (\ref%
{4.3}) holds.}

\emph{2. Let }$\lambda _{1}=\lambda _{1}\left( \Omega ,\alpha
,T,M,c,R\right) \geq 1$\emph{\ be the number of Theorem 4.3. Consider the
number }$\lambda _{2}=\lambda _{1}\left( \Omega ,\alpha ,T,M,c,3R\right)
\geq \lambda _{1}$\emph{. Then for any }$\lambda \geq \lambda _{2}$\emph{\
and for any choice of the regularization parameter }$\beta \in \left[
2e^{-\lambda \left( T/2\right) ^{1+\alpha }},1\right) $\emph{\ the
functional }$I_{\lambda ,\beta }$\emph{\ is strongly convex on the set }$%
\overline{B_{0}\left( 2R\right) }$\emph{, has unique minimizer }%
\begin{equation}
\left( \widetilde{u}_{\min ,\lambda ,\beta },\widetilde{m}_{\min ,\lambda
,\beta }\right) \in \overline{B_{0}\left( 2R\right) }  \label{4.121}
\end{equation}%
\emph{\ on this set, and the analog of (\ref{4.5}) holds.}

\emph{3. Let }$\rho \in \left( 0,1\right) $\emph{\ be an arbitrary number.
Choose the number }$\gamma =\gamma \left( \rho \right) $\emph{\ as}%
\begin{equation}
\gamma =\gamma \left( \rho \right) =\left( \frac{\rho }{3-\rho }\right)
^{1/\left( 1+\alpha \right) }\in \left( 0,\left( \frac{1}{2}\right)
^{1/\left( 1+\alpha \right) }\right) .  \label{4.13}
\end{equation}%
\emph{Choose the number }$\delta _{0}=\delta _{0}\left( \Omega ,\alpha ,\mu
,T,M,c,R,\rho \right) \in \left( 0,1\right) $\emph{\ so small that }%
\begin{equation}
\ln \left\{ \delta _{0}^{-2/\left[ \left( T/2\right) ^{1+\alpha }\left(
1+\gamma ^{1+\alpha }\right) \right] }\right\} \geq \lambda _{2}>1.
\label{4.14}
\end{equation}%
\emph{For any }$\delta \in \left( 0,\delta _{0}\right) $\emph{\ choose in }$%
I_{\lambda ,\beta }$\emph{\ parameters }$\lambda =\lambda \left( \delta
,\rho ,\alpha \right) \geq \lambda _{2}$\emph{\ and }$\beta =\beta \left(
\delta ,\rho ,\alpha \right) $\emph{\ as}%
\begin{equation}
\lambda \left( \delta ,\rho ,\alpha \right) =\ln \left\{ \delta ^{-2/\left[
\left( T/2\right) ^{1+\alpha }\left( 1+\gamma ^{1+\alpha }\right) \right]
}\right\} ,  \label{4.15}
\end{equation}

\begin{equation}
\beta =2e^{-\lambda \left( T/2\right) ^{1+\alpha }}.  \label{4.16}
\end{equation}%
\emph{Denote }%
\begin{equation}
\left( \overline{u}_{\min ,\lambda ,\beta },\overline{m}_{\min ,\lambda
,\beta }\right) =\left( \widetilde{u}_{\min ,\lambda ,\beta }+G_{1}, 
\widetilde{m}_{\min ,\lambda ,\beta }+G_{2}\right) \in \overline{B\left(
3R\right) },  \label{4.170}
\end{equation}%
\emph{see (\ref{4.120}). Then the following accuracy estimates hold:}%
\begin{equation}
\hspace{-1.5 cm} \left\Vert \left( \overline{u}_{\min ,\lambda ,\beta
}-u^{\ast },\overline{m} _{\min ,\lambda ,\beta }-m^{\ast }\right)
\right\Vert _{H^{2,1}\left( Q_{\gamma T}\right) \times H^{2,1}\left(
Q_{\gamma T}\right) }\leq C_{1}\delta ^{1-\rho },\mbox{ }\forall \delta \in
\left( 0,\delta _{0}\right) ,  \label{4.17}
\end{equation}%
\begin{equation}
\left\Vert k_{\min ,\lambda ,\beta }-k^{\ast }\right\Vert _{L_{2}\left(
\Omega \right) }\leq C_{1}\delta ^{1-\rho },\mbox{ }\forall \delta \in
\left( 0,\delta _{0}\right) .  \label{4.18}
\end{equation}%
\emph{In} \emph{(\ref{4.18}), the function }$k^{\ast }\left( x\right) $ 
\emph{is computed via (\ref{4.70}), and the function }$k_{\min ,\lambda
,\beta }\left( x\right) $\emph{\ is computed via (\ref{3.5}) and (\ref{3.6}
). }

\emph{4. Next, based on (\ref{4.90}) and (\ref{4.17}), it is reasonable to
assume that }%
\begin{equation}
\left( \overline{u}_{\min ,\lambda ,\beta },\overline{m}_{\min ,\lambda
,\beta }\right) \in \overline{B\left( R\right) }.  \label{4.180}
\end{equation}%
\emph{Then the vector function }$\left( \overline{u}_{\min ,\lambda ,\beta },%
\overline{m}_{\min ,\lambda ,\beta }\right) $ \emph{is the unique minimizer }%
$\left( u_{\min ,\lambda ,\beta },m_{\min ,\lambda ,\beta }\right) $ \emph{%
of the functional }$J_{\lambda ,\beta }$\emph{\ on the set }$\overline{%
B\left( R\right) },$\emph{\ which is found in Theorem 4.3, i.e. }%
\begin{equation}
\left( \overline{u}_{\min ,\lambda ,\beta },\overline{m}_{\min ,\lambda
,\beta }\right) =\left( u_{\min ,\lambda ,\beta },m_{\min ,\lambda ,\beta
}\right) ,  \label{4.181}
\end{equation}%
\emph{\ and, therefore, estimate (\ref{4.17}) remains valid for }$\left(
u_{\min ,\lambda ,\beta },m_{\min ,\lambda ,\beta }\right) .$

\emph{5. (Uniqueness). There exists at most one vector function }$\left(
k^{\ast },u^{\ast },m^{\ast }\right) $\emph{\ satisfying the above
conditions.}

\textbf{Remark 4.4.} \emph{If }$\rho \approx 0$\emph{\ in (\ref{4.17}), (\ref%
{4.18}), then these are almost Lipschitz stability estimates. In fact,
stability estimates, similar with ones in (\ref{4.17}), (\ref{4.18}), were
obtained in \cite{MFG6} for a different CIP for MFGS (\ref{2.5}).}

\subsection{The gradient descent method}

\label{sec:4.4}

Assume now that%
\begin{equation}
\left( u^{\ast },m^{\ast }\right) \in B^{\ast }\left( R/3-\delta \right) ,%
\mbox{ }R/3-\delta >0,  \label{4.19}
\end{equation}%
\begin{equation}
\left( \overline{u}_{\min ,\lambda ,\beta },\overline{m}_{\min ,\lambda
,\beta }\right) =\left( u_{\min ,\lambda ,\beta },m_{\min ,\lambda ,\beta
}\right) \in B\left( R/3\right) ,  \label{4.190}
\end{equation}%
where parameters $\lambda $ and $\beta $ are the ones in (\ref{4.15}) and (%
\ref{4.16}).\emph{\ }It follows from (\ref{4.17}) that assumption (\ref%
{4.190}) is reasonable, as soon as (\ref{4.19}) is true.

We construct now the gradient descent method of the minimization of the
functional $J_{\lambda ,\beta }.$ Consider an arbitrary pair of functions%
\begin{equation}
\left( u_{0},m_{0}\right) \in B\left( R/3\right) .  \label{4.20}
\end{equation}%
Let $\mu >0$ be the step size of the gradient descent method. The iterative
sequence of this method is:%
\begin{equation}
\left( u_{n},m_{n}\right) =\left( u_{n-1},m_{n-1}\right) -\mu J_{\lambda
,\beta }^{\prime }\left( u_{n-1},m_{n-1}\right) ,\mbox{ }n=1,2,...
\label{4.21}
\end{equation}%
Note that since $J_{\lambda ,\beta }^{\prime }\in H_{0}$ by Theorem 4.3,
then all pairs $\left( u_{n},m_{n}\right) $ have the same boundary
conditions, see (\ref{4.01}).

\textbf{Theorem 4.5}.\emph{\ Let }$\lambda \geq \lambda _{2}$\emph{, where }$%
\lambda _{2}$\emph{\ was chosen in Theorem 4.4. Let }$\rho ,\gamma ,\lambda $%
\emph{\ and }$\beta $\emph{\ be the parameters chosen in Theorem 4.4}$.$%
\emph{\ Assume that conditions (\ref{4.19})-(\ref{4.20}) hold. Then there
exists a number }$\mu _{0}\in \left( 0,1\right) $\emph{\ such that for any }$%
\mu \in \left( 0,\mu _{0}\right) $ \emph{there exists a number }$\theta
=\theta \left( \mu \right) \in \left( 0,1\right) $\emph{\ such that for all }%
$n\geq 1$\emph{\ } 
\begin{equation}
\hspace{-1 cm}
\left. 
\begin{array}{c}
\left( u_{n},m_{n}\right) \in B\left( R\right) ,\mbox{ } \\ 
\left\Vert u_{n}-u^{\ast }\right\Vert _{_{H^{2,1}\left( Q_{\gamma ,T}\right)
}}+\left\Vert m_{n}-m^{\ast }\right\Vert _{_{H^{2,1}\left( Q_{\gamma
,T}\right) }}+\left\Vert k_{n}-k^{\ast }\right\Vert _{L_{2}\left( \Omega
\right) }\leq \\ 
\leq C_{1}\delta ^{1-\rho }+C_{1}\theta ^{n}\left( \left\Vert u_{\min
,\lambda ,\beta }-u_{0}\right\Vert _{H^{k_{n}}\left( Q_{T}\right)
}+\left\Vert m_{\min ,\lambda ,\beta }-m_{0}\right\Vert _{H^{k_{n}}\left(
Q_{T}\right) }\right) ,%
\end{array}%
\right.  \label{4.22}
\end{equation}%
\emph{where functions }$k_{n}\left( x\right) $\emph{\ are computed via
direct analogs of (\ref{3.5}), (\ref{3.6}) with the replacement of }$u^{\ast
}$\emph{\ with }$u_{n}.$

\textbf{Proof}. We prove this theorem, assuming that Theorems 4.3 and 4.4
hold. The first line of (\ref{4.22}), existence of numbers $\mu _{0}$ and $%
\theta \left( \mu \right) $ as well as estimate 
\begin{equation}
\hspace{-1 cm} \left. 
\begin{array}{c}
\left\Vert u_{n}-u_{\min ,\lambda ,\beta }\right\Vert _{_{H^{k_{n}}\left(
Q_{T}\right) }}+\left\Vert m_{n}-m_{\min ,\lambda ,\beta }\right\Vert
_{_{H^{k_{n}}\left( Q_{T}\right) }}+\left\Vert k_{n}-k^{\ast }\right\Vert
_{L_{2}\left( \Omega \right) }\leq \\ 
\leq C_{1}\theta ^{n}\left( \left\Vert u_{\min ,\lambda ,\beta
}-u_{0}\right\Vert _{H^{k_{n}}\left( Q_{T}\right) }+\left\Vert m_{\min
,\lambda ,\beta }-m_{0}\right\Vert _{H^{k_{n}}\left( Q_{T}\right) }\right) .%
\end{array}
\right.  \label{4.23}
\end{equation}%
follow immediately from (\ref{4.19})-(\ref{4.21}) and \cite[Theorem 6]{SAR}.
Next, triangle inequality, (\ref{4.17}), (\ref{4.18}), (\ref{4.190}) and (%
\ref{4.23}) lead to the desired result. $\ \square $

\textbf{Remark 4.4}. \emph{Since }$R>0$\emph{\ is an arbitrary number and in
(\ref{4.20}) }$\left( u_{0},m_{0}\right) $\emph{\ is an arbitrary point of
the set }$B\left( R/3\right) ,$\emph{\ then Theorem 4.5 implies the }%
\underline{global } \emph{convergence of the iterative procedure (\ref{4.21}
), see Definition 1.2. Note that the gradient descent method for a
conventional functional converges only locally.}

\section{Proof of Theorem 4.1}

\label{sec:5}

Note that $1+\alpha =\left( n_{1}+n_{2}\right) /n_{2}.$ Since $n_{1}+n_{2}$
is an even number, then 
\begin{equation}
t^{1+\alpha }=\left( -t\right) ^{1+\alpha },  \label{5.1}
\end{equation}%
and, also, $t^{1+\alpha }$ makes sense for $t<0.$ We have:%
\begin{equation}
	\hspace{-1 cm}
\left. 
\begin{array}{c}
\int\limits_{0}^{d}e^{-2\lambda t^{1+\alpha }}\left(
\int\limits_{0}^{t}f\left( \tau \right) d\tau \right) ^{2}dt= \\ 
=\left( 2\lambda \left( 1+\alpha \right) \right)
^{-1}\int\limits_{0}^{d}\partial _{t}\left( -e^{-2\lambda t^{1+\alpha
}}\right) \left( 1/t^{\alpha }\right) \left( \int\limits_{0}^{t}f\left( \tau
\right) d\tau \right) ^{2}dt= \\ 
=-\left( 2\lambda \left( 1+\alpha \right) \right) ^{-1}e^{-2\lambda
d^{1+\alpha }}d^{-\alpha }\left( \int\limits_{0}^{d}f\left( \tau \right)
d\tau \right) ^{2}- \\ 
-\left( 2\lambda \left( 1+\alpha \right) \right) ^{-1}\alpha
\int\limits_{0}^{d}e^{-2\lambda t^{1+\alpha }}\left( 1/t^{\alpha +1}\right)
\left( \int\limits_{0}^{t}f\left( \tau \right) d\tau \right) ^{2}dt+ \\ 
+\left( \lambda \left( 1+\alpha \right) \right)
^{-1}\int\limits_{0}^{d}e^{-2\lambda t^{1+\alpha }}\left( 1/t^{\alpha
}\right) f\left( \int\limits_{0}^{t}f\left( \tau \right) d\tau \right) dt.%
\end{array}%
\right.  \label{5.100}
\end{equation}%
Since terms in the third and fourth lines of (\ref{5.100}) are negative,
then (\ref{5.100}) implies:%
\begin{equation}
\hspace{-1 cm}
\hspace{-1.5cm}\int\limits_{0}^{a}e^{-2\lambda t^{1+\alpha }}\left(
\int\limits_{0}^{t}f\left( \tau \right) d\tau \right) ^{2}dt\leq \frac{1}{%
\lambda \left( 1+\alpha \right) }\int\limits_{0}^{a}e^{-2\lambda t^{1+\alpha
}}\frac{1}{t^{\alpha }}f\left( \int\limits_{0}^{t}f\left( \tau \right) d\tau
\right) dt.  \label{5.2}
\end{equation}%
Hence, applying Cauchy-Schwarz inequality to the right hand side of (\ref%
{5.2}), we obtain%
\[
\int\limits_{0}^{d}e^{-2\lambda t^{1+\alpha }}\left(
\int\limits_{0}^{t}f\left( \tau \right) d\tau \right) ^{2}dt\leq 
\]%
\begin{equation}
\hspace{-1cm}\leq \frac{1}{\lambda \left( 1+\alpha \right) }\left[
\int\limits_{0}^{d}f^{2}e^{-2\lambda t^{1+\alpha }}dt\right] ^{1/2}\left[
\int\limits_{0}^{d}e^{-2\lambda t^{1+\alpha }}\frac{1}{t^{2\alpha }}\left(
\int\limits_{0}^{t}f\left( \tau \right) d\tau \right) ^{2}dt\right] ^{1/2}.
\label{5.3}
\end{equation}%
Estimate $I,$ where%
\begin{equation}
I=\int\limits_{0}^{d}e^{-2\lambda t^{1+\alpha }}\left( 1/t^{2\alpha }\right)
\left( \int\limits_{0}^{t}f\left( \tau \right) d\tau \right) ^{2}dt.
\label{5.30}
\end{equation}%
We have 
\[
\left. 
\begin{array}{c}
I=\int\limits_{0}^{d}e^{-2\lambda t^{1+\alpha }}\left( 1/t^{2\alpha }\right)
\left( \int\limits_{0}^{t}f\left( \tau \right) d\tau \right) ^{2}dt\leq \\ 
\leq \int\limits_{0}^{d}e^{-2\lambda t^{1+\alpha }}t^{1-2\alpha }\left(
\int\limits_{0}^{t}f^{2}\left( \tau \right) d\tau \right)
dt=\int\limits_{0}^{d}f^{2}\left( \tau \right) \left[ \int\limits_{\tau
}^{d}e^{-2\lambda t^{1+\alpha }}t^{1-2\alpha }dt\right] d\tau .%
\end{array}%
\right. 
\]%
Hence, 
\begin{equation}
I\leq \int\limits_{0}^{d}f^{2}\left( \tau \right) \left[ \int\limits_{\tau
}^{d}e^{-2\lambda t^{1+\alpha }}t^{1-2\alpha }dt\right] d\tau .  \label{5.4}
\end{equation}%
Estimate now the interior integral in the right hand side of (\ref{5.4}),%
\[
\left. 
\begin{array}{c}
\int\limits_{\tau }^{d}e^{-2\lambda t^{1+\alpha }}t^{1-2\alpha }dt=\left(
2\lambda \left( 1+\alpha \right) \right) ^{-1}\int\limits_{\tau
}^{d}\partial _{t}\left( -e^{-2\lambda t^{1+\alpha }}\right) t^{1-3\alpha
}dt= \\ 
=-\left( 2\lambda \left( 1+\alpha \right) \right) ^{-1}e^{-2\lambda
d^{1+\alpha }}d^{1-3\alpha }+\left( 2\lambda \left( 1+\alpha \right) \right)
^{-1}e^{-2\lambda \tau ^{1+\alpha }}\tau ^{1-3\alpha }+ \\ 
+\left( 1-3\alpha \right) \left( 2\lambda \left( 1+\alpha \right) \right)
^{-1}\int\limits_{\tau }^{d}e^{-2\lambda t^{1+\alpha }}t^{-3\alpha }dt.%
\end{array}%
\right. 
\]%
Hence,%
\begin{equation}
\hspace{-2cm}\int\limits_{\tau }^{d}e^{-2\lambda t^{1+\alpha }}t^{1-2\alpha
}dt\leq \frac{1}{2\lambda \left( 1+\alpha \right) }e^{-2\lambda \tau
^{1+\alpha }}\tau ^{1-3\alpha }+\frac{1-3\alpha }{2\lambda \left( 1+\alpha
\right) }\int\limits_{\tau }^{d}e^{-2\lambda t^{1+\alpha }}t^{-3\alpha }dt.
\label{5.5}
\end{equation}%
Next, since the function $e^{-2\lambda t^{1+\alpha }}$ is decreasing with
respect to $t$, then%
\[
\frac{1-3\alpha }{2\lambda \left( 1+\alpha \right) }\int\limits_{\tau
}^{d}e^{-2\lambda t^{1+\alpha }}t^{-3\alpha }dt\leq \frac{1-3\alpha }{%
2\lambda \left( 1+\alpha \right) }e^{-2\lambda \tau ^{1+\alpha
}}\int\limits_{\tau }^{d}t^{-3\alpha }dt= 
\]%
\[
=\frac{1}{2\lambda \left( 1+\alpha \right) }e^{-2\lambda \tau ^{1+\alpha
}}\left( d^{1-3\alpha }-\tau ^{1-3\alpha }\right) . 
\]%
Hence, in (\ref{5.5})%
\[
\hspace{-1.5cm}\int\limits_{\tau }^{d}e^{-2\lambda t^{1+\alpha
}}t^{1-2\alpha }dt\leq \frac{1}{2\lambda \left( 1+\alpha \right) }%
e^{-2\lambda \tau ^{1+\alpha }}\tau ^{1-3\alpha }+\frac{1-3\alpha }{2\lambda
\left( 1+\alpha \right) }\int\limits_{\tau }^{d}e^{-2\lambda t^{1+\alpha
}}t^{-3\alpha }dt\leq 
\]%
\[
\hspace{-2cm}\leq \frac{1}{2\lambda \left( 1+\alpha \right) }e^{-2\lambda
\tau ^{1+\alpha }}\tau ^{1-3\alpha }+\frac{1}{2\lambda \left( 1+\alpha
\right) }e^{-2\lambda \tau ^{1+\alpha }}d^{1-3\alpha }-\frac{1}{2\lambda
\left( 1+\alpha \right) }e^{-2\lambda \tau ^{1+\alpha }}\tau ^{1-3\alpha }= 
\]%
\[
=\frac{1}{2\lambda \left( 1+\alpha \right) }e^{-2\lambda \tau ^{1+\alpha
}}d^{1-3\alpha }. 
\]%
Thus, we have proven that%
\[
\int\limits_{\tau }^{d}e^{-2\lambda t^{1+\alpha }}t^{1-2\alpha }dt\leq \frac{%
d^{1-3\alpha }}{2\lambda \left( 1+\alpha \right) }e^{-2\lambda \tau
^{1+\alpha }}. 
\]%
Hence, (\ref{5.30}) and (\ref{5.4}) imply:%
\[
\hspace{-2 cm}
\left[ \int\limits_{0}^{d}e^{-2\lambda t^{1+\alpha }}\frac{1}{t^{2\alpha }}%
\left( \int\limits_{0}^{t}f\left( \tau \right) d\tau \right) ^{2}dt\right]
^{1/2}\leq \frac{d^{\left( 1-3\alpha \right) /2}}{\lambda ^{1/2}\left(
2\left( 1+\alpha \right) \right) ^{1/2}}\left(
\int\limits_{0}^{d}f^{2}\left( t\right) e^{-2\lambda t^{1+\alpha }}dt\right)
^{1/2}. 
\]%
Substituting this in (\ref{5.3}), we obtain%
\begin{equation}
\int\limits_{0}^{d}e^{-2\lambda t^{1+\alpha }}\left(
\int\limits_{0}^{t}f\left( \tau \right) d\tau \right) ^{2}dt\leq \frac{%
d^{\left( 1-3\alpha \right) /2}}{\lambda ^{3/2}\sqrt{2}\left( 1+\alpha
\right) ^{3/2}}\int\limits_{0}^{d}f^{2}e^{-2\lambda t^{1+\alpha }}dt.
\label{5.6}
\end{equation}

Next, using (\ref{5.1}) and (\ref{5.6}), we obtain 
\begin{eqnarray}
\int\limits_{-d}^{0}e^{-2\lambda t^{1+\alpha }}\left(
\int\limits_{t}^{0}f\left( \tau \right) d\tau \right)
^{2}dt=\int\limits_{0}^{d}e^{-2\lambda t^{1+\alpha }}\left(
\int\limits_{0}^{t}f\left( -\tau \right) d\tau \right) ^{2}dt\leq  \nonumber
\end{eqnarray}%
\begin{equation}
\leq \frac{1}{\lambda ^{3/2}}\cdot \frac{d^{\left( 1-3\alpha \right) /2}}{ 
\sqrt{2}\left( 1+\alpha \right) ^{3/2}}\int\limits_{0}^{d}e^{-2\lambda
t^{1+\alpha }}f^{2}\left( -t\right) dt=  \label{5.7}
\end{equation}%
\begin{eqnarray}
=\frac{1}{\lambda ^{3/2}}\cdot \frac{d^{\left( 1-3\alpha \right) /2}}{\sqrt{%
2 }\left( 1+\alpha \right) ^{3/2}}\int\limits_{-d}^{0}e^{-2\lambda
t^{1+\alpha }}f^{2}\left( t\right) dt.  \nonumber
\end{eqnarray}%
Thus, (\ref{5.7}) implies%
\begin{equation}
\hspace{-1 cm} \int\limits_{-d}^{0}e^{-2\lambda t^{1+\alpha }}\left(
\int\limits_{t}^{0}f\left( \tau \right) d\tau \right) ^{2}dt\leq \frac{1}{
\lambda ^{3/2}}\cdot \frac{d^{\left( 1-3\alpha \right) /2}}{\sqrt{2}\left(
1+\alpha \right) ^{3/2}}\int\limits_{-d}^{0}e^{-2\lambda t^{1+\alpha
}}f^{2}\left( t\right) dt.  \label{5.8}
\end{equation}%
Summing up (\ref{5.6}) and (\ref{5.8}), we obtain the target estimate (\ref%
{4.1}). \ \ $\square $

\section{Proof of Theorem 4.3}

\label{sec:6.0}

To simplify the presentation, we prove this theorem only for the case (\ref%
{2.3}). The case (\ref{2.2}) is simpler than (\ref{2.3}), see (\ref{2.4}).
It follows from (\ref{2.3}), (\ref{2.4}) and \cite[Lemma 3.2]{MFG6} that 
\begin{equation}
\hspace{-1.5cm}\int\limits_{Q_{T}}\left( \int\limits_{\Omega }K_{2}\left(
x,y\right) z\left( y,t\right) dy\right) ^{2}\varphi _{\lambda }\left(
x,t\right) dxdt\leq C_{1}\int\limits_{Q_{T}}z^{2}\varphi _{\lambda }dxdt,%
\mbox{ }\forall z\in L_{2}\left( Q_{T}\right) .  \label{6.01}
\end{equation}%
Let $\left( u_{1},m_{1}\right) ,\left( u_{2},m_{2}\right) \in \overline{%
B\left( R\right) }$ be two arbitrary pairs of functions. Denote $\left(
w,q\right) =\left( u_{2},m_{2}\right) -\left( u_{1},m_{1}\right) .$ Then
triangle inequality, (\ref{3.8}), (\ref{3.80}) and (\ref{4.10}) imply%
\begin{equation}
\left( w,q\right) \in \overline{B_{0}\left( 2R\right) }.  \label{6.02}
\end{equation}%
Also,%
\begin{equation}
J_{\lambda ,\beta }\left( u_{2},m_{2}\right) -J_{\lambda ,\beta }\left(
u_{1},m_{1}\right) =J_{\lambda ,\beta }\left( u_{1}+w,m_{1}+q\right)
-J_{\lambda ,\beta }\left( u_{1},m_{1}\right) .  \label{6.03}
\end{equation}%
Using (\ref{6.03}), we analyze now functionals $J_{1,\lambda }$ and $%
J_{2,\lambda }$ in (\ref{3.11}).

\subsection{Analysis of $J_{1,\protect\lambda }$}

\label{sec:6.01}

First, we consider the operator $L_{1}\left( u_{1}+w,m_{1}+q\right) $ in (%
\ref{3.61}) and separate its linear $L_{1,\mbox{		lin}}$ and nonlinear $L_{1,%
\mbox{nonlin}}$ parts with respect to $\left( w,q\right) .$ We drop here
dependence of $L_{1}$ on $v_{0},p_{0}$ for brevity. \ \ We have:%
\begin{equation}
\left. 
\begin{array}{c}
L_{1}\left( u_{1}+w,m_{1}+q\right) =L_{1}\left( u_{1},m_{1}\right) + \\ 
+L_{1,\mbox{lin}}\left( u_{1},m_{1},w,q\right) +L_{1,\mbox{nonlin}}\left(
u,m,w,q\right) .%
\end{array}%
\right.  \label{6.04}
\end{equation}%
Hence, 
\begin{equation}
\hspace{-2.8cm}\left. 
\begin{array}{c}
\left[ L_{1}\left( u_{1}+w,m_{1}+q\right) \right] ^{2}-\left[ L_{1}\left(
u_{1},m_{1}\right) \right] ^{2}= \\ 
=2L_{1}\left( u_{1},m_{1}\right) L_{1,\mbox{lin}}\left(
u_{1},m_{1},w,q\right) + \\ 
+\left[ L_{1,\mbox{lin}}\left( u_{1},m_{1},w,q\right) \right] ^{2}+ \\ 
+2L_{1}\left( u_{1},m_{1}\right) L_{1,\mbox{nonlin}}\left( u,m,w,q\right) +
\\ 
+2L_{1,\mbox{lin}}\left( u_{1},m_{1},w,q\right) L_{1,\mbox{nonlin}}\left(
u,m,w,q\right) + \\ 
+\left[ L_{1,\mbox{nonlin}}\left( u,m,w,q\right) \right] ^{2}.%
\end{array}%
\right.  \label{6.05}
\end{equation}%
We have: 
\begin{equation}
\hspace{-1cm}\left. 
\begin{array}{c}
L_{1,\mbox{lin}}\left( u,m,w,q\right) =w_{t}+\Delta w{-r\nabla w}\left(
\int\limits_{T/2}^{t}\nabla u_{1}(x,\tau )d\tau +\nabla v_{0}\left( x\right)
\right) - \\ 
-{r\nabla u}_{1}\int\limits_{T/2}^{t}\nabla w(x,\tau )d\tau - \\ 
-f\left( x\right) \left( w\left( x,t\right)
-\int\limits_{T/2}^{t}w_{t}\left( x,\tau \right) d\tau \right)
\int\limits_{\Omega }K\left( x,y\right) m\left( y,t\right) dy- \\ 
-f\left( x\right) \left[ \left( u\left( x,t\right)
-\int\limits_{T/2}^{t}u_{t}\left( x,\tau \right) d\tau \right) +F\left(
x\right) \right] \int\limits_{\Omega }K\left( x,y\right) q\left( y,t\right)
dy- \\ 
-sq-s_{t}\int\limits_{T/2}^{t}q(x,\tau )d\tau ,\mbox{ }\left( x,t\right) \in
Q_{T}.%
\end{array}%
\right.  \label{6.06}
\end{equation}

It is space consuming to write the explicit form of the term $L_{1,%
\mbox{		nonlin}}\left( u,m,w,q\right) .$ Nevertheless, using (\ref{2.40}), (%
\ref{3.61}), (\ref{3.7})-(\ref{3.8}), (\ref{4.10}), (\ref{6.02}), (\ref{6.04}%
)-(\ref{6.06}) and Cauchy-Schwarz inequality we obtain after routine
manipulations:%
\begin{equation}
\hspace{-1cm}\left. 
\begin{array}{c}
\left[ L_{1}\left( u_{1}+w,m_{1}+q\right) \right] ^{2}-\left[ L_{1}\left(
u_{1},m_{1}\right) \right] ^{2}- \\ 
-2L_{1}\left( u_{1},m_{1}\right) L_{1,\mbox{lin}}\left(
u_{1},m_{1},w,q\right) \geq \\ 
\geq \left( 1/2\right) \left( w_{t}+\Delta w\right) ^{2}-C_{1}\left( \nabla
w\right) ^{2}- \\ 
-C_{1}\left[ \left( \int\limits_{T/2}^{t}\left( \nabla w\right) d\tau
\right) ^{2}+\left( \int\limits_{T/2}^{t}w_{t}d\tau \right) ^{2}\right]
-C_{1}\left( \int\limits_{\Omega }K_{2}\left( x,y\right) q\left( y,t\right)
dy\right) ^{2}.%
\end{array}%
\right.  \label{6.07}
\end{equation}%
Consider the functional $J_{1,\lambda ,\mbox{lin}}\left(
u_{1},m_{1},w,q\right) $ acting on $\left( w,q\right) \in H_{0}.$ This is a
linear and bounded functional defined as:%
\begin{equation}
\hspace{-1.5 cm}
\left. 
\begin{array}{c}
J_{1,\lambda ,\mbox{lin}}\left( u_{1},m_{1},w,q\right)
=2\int\limits_{Q_{T}}L_{1}\left( u_{1},m_{1}\right) L_{1,\mbox{lin}}\left(
u_{1},m_{1},w,q\right) \varphi _{\lambda }\left( x,t\right) dxdt, \\ 
J_{1,\lambda ,\mbox{lin}}\left( u_{1},m_{1},w,q\right) :H_{0}\rightarrow 
\mathbb{R}.%
\end{array}%
.\right.  \label{6.08}
\end{equation}%
Therefore, by the Riesz theorem, there exists unique point $J_{1,\lambda ,%
\mbox{lin}}^{\prime }\left( u_{1},m_{1}\right) \in H_{0}$ such that 
\begin{equation}
J_{1,\lambda ,\mbox{lin}}\left( u_{1},m_{1},w,q\right) =\left[ J_{1,\lambda
}^{\prime }\left( u_{1},m_{1}\right) ,\left( w,q\right) \right] ,\mbox{ }%
\forall \left( w,q\right) \in H_{0},  \label{6.09}
\end{equation}%
see Remark 3.1 for $\left[ ,\right] $). Using (\ref{6.04})-(\ref{6.06}) and
considerations, which are similar to any of our above cited works on the
convexification \cite{Bak,KLpar,SAR,KL,MFG7}, we can prove that $%
J_{1,\lambda ,\mbox{lin}}^{\prime }\left( u_{1},m_{1}\right) $ is actually
the Fr\'{e}chet derivative of the functional $J_{1,\lambda ,\mbox{lin}%
}\left( u_{1},m_{1}\right) :\overline{B\left( R\right) }\rightarrow \mathbb{R%
}$ at the point $\left( u_{1},m_{1}\right) .$ The Lipschitz continuity
property (\ref{4.3}) of $J_{1,\lambda }^{\prime }\left( u,m\right) $ is
rather easy to prove, similarly with the proof of either Theorem 3.1 of \cite%
{Bak} or Theorem 5.3.1 of \cite{KL}.\ Hence, we omit the proof of (\ref{4.3}%
).

Hence, by (\ref{3.11}) and (\ref{6.07})-(\ref{6.09})%
\[
e^{-2\lambda b^{2}}\left\{ J_{1,\lambda }\left( u_{1}+w,m_{1}+q\right)
-J_{1,\lambda }\left( u_{1},m_{1}\right) -\left[ J_{1,\lambda }^{\prime
}\left( u_{1},m_{1}\right) ,\left( w,q\right) \right] \right\} \geq 
\]%
\[
\geq \left( 1/2\right) e^{-2\lambda b^{2}}\int\limits_{Q_{T}}\left(
w_{t}+\Delta w\right) ^{2}\varphi _{\lambda }dxdt-C_{1}e^{-2\lambda
b^{2}}\int\limits_{Q_{T}}\left( \nabla w\right) ^{2}\varphi _{\lambda }dxdt- 
\]%
\begin{equation}
-C_{1}e^{-2\lambda b^{2}}\int\limits_{Q_{T}}\left[ \left(
\int\limits_{T/2}^{t}\left( \nabla w\right) d\tau \right) ^{2}+\left(
\int\limits_{T/2}^{t}w_{t}d\tau \right) ^{2}\right] \varphi _{\lambda }dxdt-
\label{6.010}
\end{equation}%
\[
-C_{1}e^{-2\lambda b^{2}}\int\limits_{Q_{T}}\left( \int\limits_{\Omega
}K_{2}\left( x,y\right) q\left( y,t\right) dy\right) ^{2}\varphi _{\lambda
}dxdt. 
\]

Let $\lambda _{0}\geq 1$ be the number of Theorem 4.2. Apply now Carleman
estimate (\ref{4.2}) to the first term in the second line of (\ref{6.010}).
Also, use (\ref{6.01}). We obtain that there exists a sufficiently large
number $\lambda _{1}=\lambda _{1}\left( \Omega ,\alpha ,\mu ,T,M,c,R\right)
\geq \lambda _{0}$ such that%
\[
e^{-2\lambda b^{2}}\left\{ J_{1,\lambda }\left( u_{1}+w,m_{1}+q\right)
-J_{1,\lambda }\left( u_{1},m_{1}\right) -\left[ J_{1,\lambda }^{\prime
}\left( u_{1},m_{1}\right) ,\left( w,q\right) \right] \right\} \geq 
\]%
\[
\geq \frac{C_{1}}{\lambda }e^{-2\lambda b^{2}}\int\limits_{Q_{T}}\left(
w_{t}^{2}+\sum\limits_{i,j=1}^{n}w_{x_{i}x_{j}}^{2}\right) \varphi _{\lambda
}dxdt+ 
\]%
\begin{equation}
+C_{1}e^{-2\lambda b^{2}}\int\limits_{Q_{T}}\left( \lambda \left( \nabla
w\right) ^{2}+\lambda ^{3}w^{2}\right) \varphi _{\lambda }dxdt-
\label{6.100}
\end{equation}%
\[
-C_{1}e^{-2\lambda b^{2}}\int\limits_{Q_{T}}\left[ \left(
\int\limits_{T/2}^{t}\left( \nabla w\right) d\tau \right) ^{2}+\left(
\int\limits_{T/2}^{t}w_{t}d\tau \right) ^{2}\right] \varphi _{\lambda }dxdt- 
\]%
\[
-C_{1}e^{-2\lambda b^{2}}\int\limits_{Q_{T}}q^{2}\varphi _{\lambda }dxdt- 
\]%
\[
-C_{1}\left( \left\Vert w\left( x,T\right) \right\Vert _{H^{1}\left( \Omega
\right) }^{2}+\left\Vert w\left( x,0\right) \right\Vert _{H^{1}\left( \Omega
\right) }^{2}\right) \exp \left[ -2\lambda \left( T/2\right) ^{1+\alpha }%
\right] , 
\]%
\[
\forall \left( u_{1},m_{1}\right) ,\forall \left( u_{1}+w,m_{1}+q\right) \in 
\overline{B\left( R\right) },\mbox{ }\forall \lambda \geq \lambda _{1}. 
\]

We now need to estimate from the above the term with the Volterra integrals
in (\ref{6.100}): 
\begin{equation}
\int\limits_{Q_{T}}\left[ \left( \int\limits_{T/2}^{t}\left( \nabla w\right)
d\tau \right) ^{2}+\left( \int\limits_{T/2}^{t}w_{t}d\tau \right) ^{2}\right]
\varphi _{\lambda }dxdt.  \label{6.101}
\end{equation}%
If we would have only the term with $\nabla w,$ then the standard Carleman
Weight Function with $1+\alpha =2$ in (\ref{3.10}) would be sufficient
Indeed, in this case we would have \cite[Lemma 3.1.1]{KT}, \cite[Lemma 3.1.1]%
{KL}%
\begin{equation}
\left. 
\begin{array}{c}
\int\limits_{Q_{T}}\left( \int\limits_{T/2}^{t}\left( \nabla w\right) d\tau
\right) ^{2}\exp \left[ 2\lambda \left( x_{1}^{2}-\left( t-T/2\right)
^{2}\right) \right] dxdt\leq \\ 
\leq C_{1}\left( 1/\lambda \right) \int\limits_{Q_{T}}\left( \nabla w\right)
^{2}\exp \left[ 2\lambda \left( x_{1}^{2}-\left( t-T/2\right) ^{2}\right) %
\right] dxdt,%
\end{array}%
\right.  \label{6.102}
\end{equation}%
also see Remark 4.1. Next, since we have $+\lambda \left( \nabla w\right)
^{2}$ in (\ref{6.100}) and since $\lambda >>1/\lambda $ for sufficiently
large $\lambda ,$ then the term in the second line of (\ref{6.102}) would be
dominated in (\ref{6.100}) by the term with $+\lambda \left( \nabla w\right)
^{2}.$ Thus, we are especially concerned with the second term in (\ref{6.101}%
) since the term with $w_{t}^{2}$ is multiplied by the multiplier $1/\lambda 
$ in (\ref{6.100}), which tends to zero as $\lambda \rightarrow \infty .$ In
fact, we need Theorem 4.1 exactly due to our concern with $w_{t}$.

Using (\ref{3.10}) and Theorem 4.1, we obtain 
\[
\hspace{-1.5cm}\int\limits_{Q_{T}}\left[ \left( \int\limits_{T/2}^{t}\left(
\nabla w\right) d\tau \right) ^{2}+\left( \int\limits_{T/2}^{t}w_{t}d\tau
\right) ^{2}\right] \varphi _{\lambda }dxdt\leq \frac{C_{1}}{\lambda ^{3/2}}%
\int\limits_{Q_{T}}\left[ \left( \nabla w\right) ^{2}+w_{t}^{2}\right]
\varphi_{\lambda }dxdt. 
\]%
Hence, since 
\[
\frac{1}{\lambda ^{3/2}}<<\frac{1}{\lambda }\mbox{ and }\frac{1}{\lambda
^{3/2}}<<\lambda ,\mbox{ }\forall \lambda \geq \lambda _{1}, 
\]%
then (\ref{6.100}) becomes:%
\begin{equation}
\hspace{-2cm}\left. 
\begin{array}{c}
e^{-2\lambda b^{2}}S\left( J_{1,\lambda }\right) = \\ 
=e^{-2\lambda b^{2}}\left\{ J_{1,\lambda }\left( u_{1}+w,m_{1}+q\right)
-J_{1,\lambda }\left( u_{1},m_{1}\right) -\left[ J_{1,\lambda }^{\prime
}\left( u_{1},m_{1}\right) ,\left( w,q\right) \right] \right\} \geq \\ 
\geq \left( C_{1}/\lambda \right) e^{-2\lambda
b^{2}}\int\limits_{Q_{T}}\left(
w_{t}^{2}+\sum\limits_{i,j=1}^{n}w_{x_{i}x_{j}}^{2}\right) \varphi _{\lambda
}dxdt+ \\ 
+C_{1}e^{-2\lambda b^{2}}\int\limits_{Q_{T}}\left( \lambda \left( \nabla
w\right) ^{2}+\lambda ^{3}w^{2}\right) \varphi _{\lambda }dxdt- \\ 
-C_{1}e^{-2\lambda b^{2}}\int\limits_{Q_{T}}q^{2}\varphi _{\lambda }dxdt- \\ 
-C_{1}\left( \left\Vert w\left( x,T\right) \right\Vert _{H^{1}\left( \Omega
\right) }^{2}+\left\Vert w\left( x,0\right) \right\Vert _{H^{1}\left( \Omega
\right) }^{2}\right) \exp \left[ -2\lambda \left( T/2\right) ^{1+\alpha }%
\right] , \\ 
\forall \left( u_{1},m_{1}\right) ,\forall \left( u_{1}+w,m_{1}+q\right) \in 
\overline{B\left( R\right) },\mbox{ }\forall \lambda \geq \lambda _{1}.%
\end{array}%
\right.  \label{6.011}
\end{equation}

\subsection{Analysis of $\protect\lambda ^{3/2}J_{1,\protect\lambda }+J_{2, 
\protect\lambda }$ and of $J_{\protect\lambda ,\protect\beta }$}

\label{sec:6.02}

Using formula (\ref{3.62}) and considerations, which are completely similar
with ones of the previous subsection, we obtain the following analog of (\ref%
{6.011}):%
\begin{equation}
\hspace{-2cm}\left. 
\begin{array}{c}
e^{-2\lambda b^{2}}S\left( J_{2,\lambda }\right) = \\ 
=e^{-2\lambda b^{2}}\left\{ J_{2,\lambda }\left( u_{1}+w,m_{1}+q\right)
-J_{2,\lambda }\left( u_{1},m_{1}\right) -\left[ J_{2,\lambda }^{\prime
}\left( u_{1},m_{1}\right) ,\left( w,q\right) \right] \right\} \geq \\ 
\geq \left( C_{1}/\lambda \right) e^{-2\lambda
b^{2}}\int\limits_{Q_{T}}\left(
q_{t}^{2}+\sum\limits_{i,j=1}^{n}q_{x_{i}x_{j}}^{2}\right) \varphi _{\lambda
}dxdt+ \\ 
+C_{1}e^{-2\lambda b^{2}}\int\limits_{Q_{T}}\left( \lambda \left( \nabla
q\right) ^{2}+\lambda ^{3}q^{2}\right) \varphi _{\lambda }dxdt- \\ 
-C_{1}e^{-2\lambda b^{2}}\int\limits_{Q_{T}}\left( \Delta w\right)
^{2}\varphi _{\lambda }dxdt-e^{-2\lambda b^{2}}\int\limits_{Q_{T}}\left(
\nabla w\right) ^{2}\varphi _{\lambda }dxdt- \\ 
-C_{1}\left( \left\Vert q\left( x,T\right) \right\Vert _{H^{1}\left( \Omega
\right) }^{2}+\left\Vert q\left( x,0\right) \right\Vert _{H^{1}\left( \Omega
\right) }^{2}\right) \exp \left[ -2\lambda \mu \left( T/2\right) ^{1+\alpha }%
\right] , \\ 
\forall \left( u_{1},m_{1}\right) ,\forall \left( u_{1}+w,m_{1}+q\right) \in 
\overline{B\left( R\right) },\mbox{ }\forall \lambda \geq \lambda _{1}.%
\end{array}%
\right.  \label{6.012}
\end{equation}%
Multiply both sides of (\ref{6.011}) by $\lambda ^{3/2}$ and sum up with (%
\ref{6.012}). We obtain%
\begin{equation}
\hspace{-1cm}\left. 
\begin{array}{c}
e^{-2\lambda b^{2}}\left[ \lambda ^{3/2}S\left( J_{1,\lambda }\right)
+S\left( J_{2,\lambda }\right) \right] \geq \\ 
\geq C_{1}e^{-2\lambda b^{2}}\sqrt{\lambda }\int\limits_{Q_{T}}\left(
w_{t}^{2}+\sum\limits_{i,j=1}^{n}w_{x_{i}x_{j}}^{2}\right) \varphi _{\lambda
}dxdt+ \\ 
+e^{-2\lambda b^{2}}\left( C_{1}/\lambda \right) \int\limits_{Q_{T}}\left(
q_{t}^{2}+\sum\limits_{i,j=1}^{n}q_{x_{i}x_{j}}^{2}\right) \varphi _{\lambda
}dxdt+ \\ 
+C_{1}e^{-2\lambda b^{2}}\int\limits_{Q_{T}}\left[ \lambda ^{5/2}\left(
\nabla w\right) ^{2}+\lambda ^{9/2}w^{2}+\lambda \left( \nabla q\right)
^{2}+\lambda ^{3}q^{2}\right] \varphi _{\lambda }dxdt- \\ 
-C_{1}\left( \left\Vert w\left( x,T\right) \right\Vert _{H^{1}\left( \Omega
\right) }^{2}+\left\Vert w\left( x,0\right) \right\Vert _{H^{1}\left( \Omega
\right) }^{2}\right) \lambda ^{3/2}\exp \left[ -2\lambda \left( T/2\right)
^{1+\alpha }\right] - \\ 
-C_{1}\left( \left\Vert q\left( x,T\right) \right\Vert _{H^{1}\left( \Omega
\right) }^{2}+\left\Vert q\left( x,0\right) \right\Vert _{H^{1}\left( \Omega
\right) }^{2}\right) \lambda ^{3/2}\exp \left[ -2\lambda \left( T/2\right)
^{1+\alpha }\right] , \\ 
\forall \left( u_{1},m_{1}\right) ,\forall \left( u_{1}+w,m_{1}+q\right) \in 
\overline{B\left( R\right) },\mbox{ }\forall \lambda \geq \lambda _{1}.%
\end{array}%
\right.  \label{6.013}
\end{equation}%
Consider now an arbitrary number $\gamma \in \left( 0,1\right) .$ Since by
the last line of (\ref{2.1}) $Q_{\gamma T}\subset Q_{T},$ then, making
inequality (\ref{6.013}) stronger and using (\ref{3.7}) and (\ref{4.01}), we
obtain%
\begin{equation}
\left. 
\begin{array}{c}
\lambda ^{3/2}\exp \left[ -2\lambda \left( T/2\right) ^{1+\alpha }\right]
\left\Vert \left( w,q\right) \right\Vert _{H}^{2}+ \\ 
+\left( \lambda ^{3/2}S\left( J_{1,\lambda }\right) +S\left( J_{2,\lambda
}\right) \right) \geq \\ 
\geq C_{1}\left[ \exp \left( -2\lambda \left( \gamma T/2\right) ^{1+\alpha
}\right) /\lambda \right] \left\Vert \left( w,q\right) \right\Vert
_{H^{2,1}\left( Q_{\gamma T}\right) \times H^{2,1}\left( Q_{\gamma T}\right)
}^{2}, \\ 
\forall \left( u_{1},m_{1}\right) ,\forall \left( u_{1}+w,m_{1}+q\right) \in 
\overline{B\left( R\right) },\mbox{ }\forall \lambda \geq \lambda _{1},%
\mbox{
		}\forall \gamma \in \left( 0,1\right) .%
\end{array}%
\right.  \label{6.014}
\end{equation}%
Consider now the functional $J_{\lambda ,\beta }\left( u,m\right) $ in (\ref%
{3.11}). By (\ref{6.014})%
\begin{equation}
\hspace{-1cm}\left. 
\begin{array}{c}
\lambda ^{3/2}\exp \left[ -2\lambda \left( T/2\right) ^{1+\alpha }\right]
\left\Vert \left( w,q\right) \right\Vert _{H}^{2}+ \\ 
+J_{\lambda ,\beta }\left( u_{1}+w,m_{1}+q\right) -J_{\lambda ,\beta }\left(
u_{1},m_{1}\right) -\left[ J_{\lambda ,\beta }^{\prime }\left(
u_{1},m_{1}\right) ,\left( w,q\right) \right] \geq \\ 
\geq C_{1}\left\{ \exp \left[ -2\lambda \left( \gamma T/2\right) ^{1+\alpha }%
\right] /\lambda \right\} \left\Vert \left( w,q\right) \right\Vert
_{H^{2,1}\left( Q_{\gamma T}\right) \times H^{2,1}\left( Q_{\gamma T}\right)
}^{2}+ \\ 
+\beta \left\Vert \left( w,q\right) \right\Vert _{H}^{2}, \\ 
\forall \left( u_{1},m_{1}\right) ,\left( u_{1}+w,m_{1}+q\right) \in 
\overline{B\left( R\right) },\mbox{ }\forall \lambda \geq \lambda _{1},%
\mbox{
		}\forall \gamma \in \left( 0,1\right) .%
\end{array}%
\right.  \label{6.015}
\end{equation}%
Next, since $\beta \in \left[ 2\exp \left[ -\lambda \left( T/2\right)
^{1+\alpha }\right] ,1\right) ,$ then $\beta /2>\lambda ^{3/2}\exp \left[
-2\lambda \left( T/2\right) ^{1+\alpha }\right] $ for $\lambda \geq \lambda
_{1}.$ Hence, (\ref{6.015}) implies%
\[
\left. 
\begin{array}{c}
J_{\lambda ,\beta }\left( u_{1}+w,m_{1}+q\right) -J_{\lambda ,\beta }\left(
u_{1},m_{1}\right) -\left[ J_{\lambda ,\beta }^{\prime }\left(
u_{1},m_{1}\right) ,\left( w,q\right) \right] \geq \\ 
\begin{array}{c}
\geq C_{1}\left\{ \exp \left[ -2\lambda \left( \gamma T/2\right) ^{1+\alpha }%
\right] /\lambda \right\} \left\Vert \left( w,q\right) \right\Vert
_{H^{2,1}\left( Q_{\gamma T}\right) \times H^{2,1}\left( Q_{\gamma T}\right)
}^{2}+ \\ 
+\left( \beta /2\right) \left\Vert \left( w,q\right) \right\Vert _{H}^{2},
\\ 
\forall \left( u_{1},m_{1}\right) ,\forall \left( u_{1}+w,m_{1}+q\right) \in 
\overline{B\left( R\right) },\mbox{ }\forall \lambda \geq \lambda _{1},%
\end{array}%
\end{array}%
\right. 
\]%
which proves the strong convexity property (\ref{4.4}).

As soon as (\ref{4.4}) is established, the existence and uniqueness of the
minimizer $\left( u_{\min ,\lambda ,\beta },m_{\min ,\lambda ,\beta }\right) 
$\ of the functional $J_{\lambda ,\beta }\left( u,m\right) $ on the set $%
\overline{B\left( R\right) }$\ as well as inequality (\ref{4.5}) follow
immediately from a combination of Lemma 2.1 with Theorem 2.1 of \cite{Bak}
as well as from a combination of Lemma 5.2.1 with Theorem 5.2.1 of \cite{KL}%
. $\square $

\section{Proof of Theorem 4.4}

\label{sec:7}

It follows from (\ref{4.12}) and (\ref{4.120}) that Theorem 4.3 remains
valid for the functional $I_{\lambda ,\beta }$ for all $\lambda \geq \lambda
_{2}.$ For these values of $\lambda ,$ let $\left( \widetilde{u}_{\min
,\lambda ,\beta },\widetilde{m}_{\min ,\lambda ,\beta }\right) \in \overline{%
B_{0}\left( 2R\right) }$ be the unique minimizer of the functional $%
I_{\lambda ,\beta }$, which was found in Theorem 4.3. Hence, using (\ref{4.4}%
), we obtain%
\begin{equation}
\left. 
\begin{array}{c}
I_{\lambda ,\beta }\left( \widetilde{u}^{\ast },\widetilde{m}^{\ast }\right)
-I_{\lambda ,\beta }\left( \widetilde{u}_{\min ,\lambda ,\beta },\widetilde{m%
}_{\min ,\lambda ,\beta }\right) - \\ 
-\left[ I_{\lambda ,\beta }^{\prime }\left( \widetilde{u}_{\min ,\lambda
,\beta },\widetilde{m}_{\min ,\lambda ,\beta }\right) ,\left( \widetilde{u}%
^{\ast }-\widetilde{u}_{\min ,\lambda ,\beta },\widetilde{m}^{\ast }-%
\widetilde{m}_{\min ,\lambda ,\beta }\right) \right] \geq \\ 
\geq C_{1}\exp \left[ -2\lambda \left( \gamma T/2\right) ^{1+\alpha }\right]
\times \\ 
\times \left\Vert \left( \widetilde{u}^{\ast }-\widetilde{u}_{\min ,\lambda
,\beta },\widetilde{m}^{\ast }-\widetilde{m}_{\min ,\lambda ,\beta }\right)
\right\Vert _{H^{2,1}\left( Q_{\gamma T}\right) \times H^{2,1}\left(
Q_{\gamma T}\right) }^{2}, \\ 
\forall \lambda \geq \lambda _{2},\mbox{ }\forall \gamma \in \left(
0,1\right) .%
\end{array}%
\right.  \label{7.1}
\end{equation}%
By (\ref{4.5})%
\[
\left. 
\begin{array}{c}
-I_{\lambda ,\beta }\left( \widetilde{u}_{\min ,\lambda ,\beta },\widetilde{m%
}_{\min ,\lambda ,\beta }\right) - \\ 
-\left[ I_{\lambda ,\beta }^{\prime }\left( \widetilde{u}_{\min ,\lambda
,\beta },\widetilde{m}_{\min ,\lambda ,\beta }\right) ,\left( \widetilde{u}%
^{\ast }-\widetilde{u}_{\min ,\lambda ,\beta },\widetilde{m}^{\ast }-%
\widetilde{m}_{\min ,\lambda ,\beta }\right) \right] \leq 0.%
\end{array}%
\right. 
\]%
Hence, (\ref{7.1}) implies%
\begin{equation}
\left. 
\begin{array}{c}
C_{1}\exp \left[ 2\lambda \left( \gamma T/2\right) ^{1+\alpha }\right]
I_{\lambda ,\beta }\left( \widetilde{u}^{\ast },\widetilde{m}^{\ast }\right)
\geq \\ 
\geq \left\Vert \left( \widetilde{u}^{\ast }-\widetilde{u}_{\min ,\lambda
,\beta },\widetilde{m}^{\ast }-\widetilde{m}_{\min ,\lambda ,\beta }\right)
\right\Vert _{H^{2,1}\left( Q_{\gamma T}\right) \times H^{2,1}\left(
Q_{\gamma T}\right) }^{2}, \\ 
\forall \lambda \geq \lambda _{2},\mbox{ }\forall \lambda \geq \lambda _{2},%
\mbox{ }\forall \gamma \in \left( 0,1\right) .%
\end{array}%
\right.  \label{7.2}
\end{equation}%
Now, by (\ref{4.9})%
\begin{equation}
\hspace{-1.5cm}\left( \widetilde{u}^{\ast }-\widetilde{u}_{\min ,\lambda
,\beta },\widetilde{m}^{\ast }-\widetilde{m}_{\min ,\lambda ,\beta }\right)
=\left( u^{\ast },m^{\ast }\right) -G^{\ast }-\left[ \left( \widetilde{u}%
_{\min ,\lambda ,\beta },\widetilde{m}_{\min ,\lambda ,\beta }\right) +G-G%
\right] .  \label{7.3}
\end{equation}%
Recall that by (\ref{4.120}) 
\[
\left( \widetilde{u}_{\min ,\lambda ,\beta },\widetilde{m}_{\min ,\lambda
,\beta }\right) +G=\left( \overline{u}_{\min ,\lambda ,\beta },\overline{m}%
_{\min ,\lambda ,\beta }\right) \in \overline{B\left( 3R\right) }. 
\]%
Hence, (\ref{7.3}) implies 
\begin{equation}
\hspace{-1.5cm}\left( \widetilde{u}^{\ast }-\widetilde{u}_{\min ,\lambda
,\beta },\widetilde{m}^{\ast }-\widetilde{m}_{\min ,\lambda ,\beta }\right)
=\left( u^{\ast }-\overline{u}_{\min ,\lambda ,\beta },m^{\ast }-\overline{m}%
_{\min ,\lambda ,\beta }\right) +\left( G-G^{\ast }\right) .  \label{7.4}
\end{equation}%
On the other hand, by (\ref{4.8}) 
\begin{equation}
\left\Vert G-G^{\ast }\right\Vert _{H}\leq \delta .  \label{7.5}
\end{equation}%
Hence, (\ref{7.4}), (\ref{7.5}) and Cauchy-Schwarz inequality imply%
\[
\left. 
\begin{array}{c}
\left\Vert \left( \widetilde{u}^{\ast }-\widetilde{u}_{\min ,\lambda ,\beta
},\widetilde{m}^{\ast }-\widetilde{m}_{\min ,\lambda ,\beta }\right)
\right\Vert _{H^{2,1}\left( Q_{\gamma T}\right) \times H^{2,1}\left(
Q_{\gamma T}\right) }^{2}\geq \\ 
\geq \left( 1/2\right) \left\Vert \left( u^{\ast }-\overline{u}_{\min
,\lambda ,\beta },m^{\ast }-\overline{m}_{\min ,\lambda ,\beta }\right)
\right\Vert _{H^{2,1}\left( Q_{\gamma T}\right) \times H^{2,1}\left(
Q_{\gamma T}\right) }^{2}-\delta ^{2}.%
\end{array}%
\right. 
\]%
Hence, using (\ref{7.2}), we obtain%
\begin{equation}
\left. 
\begin{array}{c}
C_{1}\exp \left[ 2\lambda \left( \gamma T/2\right) ^{1+\alpha }\right]
I_{\lambda ,\beta }\left( \widetilde{u}^{\ast },\widetilde{m}^{\ast }\right)
+C_{1}\delta ^{2}\geq \\ 
\geq \left\Vert \left( u^{\ast }-\overline{u}_{\min ,\lambda ,\beta
},m^{\ast }-\overline{m}_{\min ,\lambda ,\beta }\right) \right\Vert
_{H^{2,1}\left( Q_{\gamma T}\right) \times H^{2,1}\left( Q_{\gamma T}\right)
}^{2}, \\ 
\forall \lambda \geq \lambda _{2},\mbox{ }\forall \gamma \in \left(
0,1\right) .%
\end{array}%
\right.  \label{7.6}
\end{equation}

Consider now $I_{\lambda ,\beta }\left( \widetilde{u}^{\ast },\widetilde{m}%
^{\ast }\right) .$ By (\ref{4.9}) and (\ref{4.12})%
\begin{equation}
\left. 
\begin{array}{c}
I_{\lambda ,\beta }\left( \widetilde{u}^{\ast },\widetilde{m}^{\ast }\right)
=J_{\lambda ,\beta }\left( \widetilde{u}^{\ast }+G_{1},\widetilde{m}^{\ast
}+G_{2}\right) = \\ 
=J_{\lambda ,\beta }\left( u^{\ast }+\left( G_{1}-G_{1}^{\ast }\right)
,m^{\ast }+\left( G_{2}-G_{2}^{\ast }\right) \right) .%
\end{array}
\right.  \label{7.7}
\end{equation}%
Next, using (\ref{3.11}), we obtain%
\begin{equation}
\left. 
\begin{array}{c}
J_{\lambda ,\beta }\left( u^{\ast }+\left( G_{1}-G_{1}^{\ast }\right)
,m^{\ast }+\left( G_{2}-G_{2}^{\ast }\right) \right) = \\ 
=e^{-2\lambda b^{2}}\lambda ^{3/2}\int\limits_{Q_{T}}\left( L_{1}\left(
u^{\ast }+\left( G_{1}-G_{1}^{\ast }\right) ,m^{\ast }+\left(
G_{2}-G_{2}^{\ast }\right) \right) \right) ^{2}\varphi _{\lambda }dxdt+ \\ 
+e^{-2\lambda b^{2}}\int\limits_{Q_{T}}\left( L_{2}\left( u^{\ast }+\left(
G_{1}-G_{1}^{\ast }\right) ,m^{\ast }+\left( G_{2}-G_{2}^{\ast }\right)
\right) \right) ^{2}\varphi _{\lambda }dxdt+ \\ 
+\beta \left\Vert \left( u^{\ast }+\left( G_{1}-G_{1}^{\ast }\right)
,m^{\ast }+\left( G_{2}-G_{2}^{\ast }\right) \right) \right\Vert _{H}^{2}.%
\end{array}
\right.  \label{7.8}
\end{equation}%
It follows from (\ref{3.100}), (\ref{4.72}), (\ref{4.8}) and (\ref{4.11})
that%
\begin{equation}
\left. 
\begin{array}{c}
e^{-2\lambda b^{2}}\lambda ^{3/2}\int\limits_{Q_{T}}\left( L_{1}\left(
u^{\ast }+\left( G_{1}-G_{1}^{\ast }\right) ,m^{\ast }+\left(
G_{2}-G_{2}^{\ast }\right) \right) \right) ^{2}\varphi _{\lambda }dxdt+ \\ 
+e^{-2\lambda b^{2}}\int\limits_{Q_{T}}\left( L_{2}\left( u^{\ast }+\left(
G_{1}-G_{1}^{\ast }\right) ,m^{\ast }+\left( G_{2}-G_{2}^{\ast }\right)
\right) \right) ^{2}\varphi _{\lambda }dxdt+ \\ 
+\beta \left\Vert \left( u^{\ast }+\left( G_{1}-G_{1}^{\ast }\right)
,m^{\ast }+\left( G_{2}-G_{2}^{\ast }\right) \right) \right\Vert _{H}^{2}\leq
\\ 
\leq C_{1}\lambda ^{3/2}\delta ^{2}+C_{1}\beta .%
\end{array}
\right.  \label{7.9}
\end{equation}%
Hence, by (\ref{7.7})-(\ref{7.9})%
\[
I_{\lambda ,\beta }\left( \widetilde{u}^{\ast },\widetilde{m}^{\ast }\right)
\leq C_{1}\lambda ^{3/2}\delta ^{2}+C_{1}\beta . 
\]%
Hence, (\ref{7.6}) implies%
\begin{equation}
\left. 
\begin{array}{c}
C_{1}\exp \left[ 3\lambda \left( \gamma T/2\right) ^{1+\alpha }\right]
\delta ^{2}+C_{1}\beta \exp \left[ 2\lambda \left( \gamma T/2\right)
^{1+\alpha }\right] \geq \\ 
\geq \left\Vert \left( u^{\ast }-\overline{u}_{\min ,\lambda ,\beta
},m^{\ast }-\overline{m}_{\min ,\lambda ,\beta }\right) \right\Vert
_{H^{2,1}\left( Q_{\gamma T}\right) \times H^{2,1}\left( Q_{\gamma T}\right)
}^{2}, \\ 
\forall \lambda \geq \lambda _{2},\mbox{ }\forall \gamma \in \left(
0,1\right) .%
\end{array}
\right.  \label{7.10}
\end{equation}%
Recall that by (\ref{4.16}) $\beta =2e^{-\lambda \left( T/2\right)
^{1+\alpha }}.$ Hence, using (\ref{7.10}), we obtain%
\begin{equation}
\left. 
\begin{array}{c}
\left\Vert \left( u^{\ast }-\overline{u}_{\min ,\lambda ,\beta },m^{\ast }- 
\overline{m}_{\min ,\lambda ,\beta }\right) \right\Vert _{H^{2,1}\left(
Q_{\gamma T}\right) \times H^{2,1}\left( Q_{\gamma T}\right) }^{2}\leq \\ 
\leq C_{1}\exp \left[ 3\lambda \left( \gamma T/2\right) ^{1+\alpha }\right]
\delta ^{2}+C_{1}\exp \left[ -\lambda \left( T/2\right) ^{1+\alpha }\left(
1-2\gamma ^{1+\alpha }\right) \right] .%
\end{array}
\right.  \label{7.11}
\end{equation}

It follows from (\ref{4.13}) that 
\begin{equation}
2\gamma ^{1+\alpha }<1.  \label{7.12}
\end{equation}%
Find $\lambda =\lambda \left( \delta ,\rho ,\alpha \right) $ such that%
\begin{equation}
\exp \left[ 3\lambda \gamma ^{1+\alpha }\left( T/2\right) ^{1+\alpha }\right]
\delta ^{2}=\exp \left[ -\lambda \left( T/2\right) ^{1+\alpha }\left(
1-2\gamma ^{1+\alpha }\right) \right] .  \label{7.13}
\end{equation}%
Hence, 
\begin{equation}
3\lambda \gamma ^{1+\alpha }\left( T/2\right) ^{1+\alpha }+\lambda \left(
T/2\right) ^{1+\alpha }\left( 1-2\gamma ^{1+\alpha }\right) =\ln \delta
^{-2}.  \label{7.14}
\end{equation}%
Hence, $\lambda =\lambda \left( \delta ,\rho ,\alpha \right) $ is as in (\ref%
{4.15}). In particular, in order to ensure that $\lambda =\lambda \left(
\delta ,\rho ,\alpha \right) $ $\geq \lambda _{2},$ we should take the
number $\delta _{0}=\delta _{0}\left( \Omega ,\alpha ,\mu ,T,M,c,R,\gamma
\right) \in \left( 0,1\right) $ so small that (\ref{4.14}) holds, and we
should also take $\delta \in \left( 0,\delta _{0}\right) .$

Using (\ref{7.11})-(\ref{7.14}) and (\ref{4.15}), we obtain 
\begin{equation}
\hspace{-1.5cm}\left\Vert \left( u^{\ast }-\overline{u}_{\min ,\lambda
,\beta },m^{\ast }-\overline{m}_{\min ,\lambda ,\beta }\right) \right\Vert
_{H^{2,1}\left( Q_{\gamma T}\right) \times H^{2,1}\left( Q_{\gamma T}\right)
}\leq C_{1}\delta ^{\left( 1-2\gamma ^{1+\alpha }\right) /\left( 1+\gamma
^{1+\alpha }\right) }.  \label{7.15}
\end{equation}%
Note that by (\ref{7.12})%
\[
\frac{1-2\gamma ^{1+\alpha }}{1+\gamma ^{1+\alpha }}\in \left( 0,1\right) . 
\]%
Hence, for any $\rho \in \left( 0,1\right) $ we can choose $\gamma =\gamma
\left( \rho \right) $ as in (\ref{4.13}). Substituting this $\gamma \left(
\rho \right) $ in (\ref{7.15}), we obtain the first target estimate (\ref%
{4.17}). Next, using (\ref{3.5}), (\ref{3.6}), (\ref{4.70}) and (\ref{4.17}%
), we obtain the second target estimate (\ref{4.18}).

Assume now that (\ref{4.180}) holds and prove (\ref{4.181}). Below in this
section 7 $\lambda =\lambda \left( \delta ,\rho ,\alpha \right) $ and as in (%
\ref{4.15}) and $\beta $ is as in (\ref{4.16}).

By (\ref{4.12}), (\ref{4.121}) and (\ref{4.170}) 
\begin{equation}
\left. 
\begin{array}{c}
J_{\lambda ,\beta }\left( \widetilde{u}_{\min ,\lambda ,\beta }+G_{1}, 
\widetilde{m}_{\min ,\lambda ,\beta }+G_{2}\right) =J_{\lambda ,\beta
}\left( \overline{u}_{\min ,\lambda ,\beta },\overline{m}_{\min ,\lambda
,\beta }\right) \leq \\ 
\leq J_{\lambda ,\beta }\left( \widetilde{u}+G_{1},\widetilde{m}
+G_{2}\right) ,\mbox{ }\forall \left( \widetilde{u},\widetilde{m}\right) \in 
\overline{B_{0}\left( 2R\right) }.%
\end{array}
\right.  \label{7.16}
\end{equation}%
Consider now an arbitrary point $\left( u_{1},m_{1}\right) \in \overline{
B\left( R\right) }.$ Then by (\ref{4.90}) and triangle inequality 
\begin{equation}
\left( \widetilde{u}_{1},\widetilde{m}_{1}\right) =\left(
u_{1}-G_{1},m_{1}-G_{2}\right) \in \overline{B_{0}\left( 2R\right) }.
\label{7.17}
\end{equation}%
Hence, using (\ref{7.16}) and (\ref{7.17}), we obtain%
\[
J_{\lambda ,\beta }\left( \overline{u}_{\min ,\lambda ,\beta },\overline{m}
_{\min ,\lambda ,\beta }\right) \leq J_{\lambda ,\beta }\left( \widetilde{u}
_{1}+G_{1},\widetilde{m}_{1}+G_{2}\right) =J_{\lambda ,\beta }\left(
u_{1},m_{1}\right) . 
\]%
Hence, (\ref{4.180}) and Theorem 4.3 imply that $\left( \overline{u}_{\min
,\lambda ,\beta },\overline{m}_{\min ,\lambda ,\beta }\right) $ is the
unique minimizer of the functional $J_{\lambda ,\beta }\left( u,m\right) $
on the set $\overline{B\left( R\right) },$ which proves (\ref{4.181}).

We now prove uniqueness of our CIP. Set $\delta =0.$ Also, for an arbitrary
number $\rho \in \left( 0,1\right) ,$ set $\gamma =\gamma \left( \rho
\right) $ as in (\ref{4.13}). Then (\ref{4.17}) and (\ref{4.18}) imply $%
\left( k^{\ast },u^{\ast },m^{\ast }\right) =\left( k_{\min ,\lambda ,\beta
},\overline{u}_{\min ,\lambda ,\beta },\overline{m}_{\min ,\lambda ,\beta
}\right) $ in the set $\Omega \times Q_{\gamma T}\times Q_{\gamma T}.$ Since
the triple $\left( k_{\min ,\lambda ,\beta },\overline{u}_{\min ,\lambda
,\beta },\overline{m}_{\min ,\lambda ,\beta }\right) $ is unique in this
set, then the triple $\left( k^{\ast },u^{\ast },m^{\ast }\right) $ is also
unique in this set. In particular, since the function $k^{\ast }=k^{\ast
}\left( x\right) ,x\in \Omega ,$ then this function is found uniquely,
independently on $\gamma .$ Given the latter, uniqueness of the vector
function $\left( u^{\ast },m^{\ast }\right) $ in the entire domain $%
Q_{T}\times Q_{T}$ follows immediately from Theorems 5 and 6 of \cite{MFG4}. 
$\square $

\section{Numerical Studies}

\label{sec:8}

In this section we describe our numerical studies of the Minimization
Problem formulated in section 3. As it is always done in numerical studies
of Ill-Posed and Inverse Problems, we need to figure out how to numerically
generate the data for our problem. More precisely, we need to numerically
generate the observation data (\ref{2.400}) of our CIP. On the first step,
we take with the coefficient $k(x)$ of our choice and numerically generate
the pair of functions $(v(x,t),p(x,t))$, which solves MFGS (\ref{2.5}). This
vector function $(v(x,t),p(x,t)),$ generates the observation data (\ref%
{2.400}). On the second step, we \textquotedblleft pretend" that we do not
know the coefficient $k(x)$ and solve the Minimization Problem with the
observation data (\ref{2.400}). Next, we compare the resulting computed
solution $k_{\mbox{comp}}(x)$ with the $k(x)$ given in the first step.

Our new procedure for data generation for MFGS (\ref{2.5}) is described in
subsection 8.1, and the numerical tests of computing the function $k_{%
\mbox{comp}}(x)$ are described in subsection 8.2.

\subsection{Numerical data generation}

\label{sec:8.1}

Choose a sufficiently smooth function $v^{\ast \ast }(x,t)$. Then we solve
the following modification of the second equation in (\ref{2.5}): 
\begin{equation}
\hspace{-0.5cm}p_{t}^{\ast \ast }(x,t)-\Delta p^{\ast \ast }(x,t){-\mbox{div}%
(r(x)p^{\ast \ast }(x,t)\nabla v^{\ast \ast }(x,t))}=0,\quad \left(
x,t\right) \in Q_{T},  \label{8.1}
\end{equation}%
with the initial condition 
\begin{equation}
p^{\ast \ast }(x,0)=p_{0}^{\ast \ast }(x),\quad x\in \Omega ,  \label{8.2}
\end{equation}%
and the Dirichlet boundary condition 
\begin{equation}
p^{\ast \ast }(x,t)\mid _{S_{T}}=g_{0,2}^{\ast \ast }(x,t),\quad \left(
x,t\right) \in S_{T},  \label{8.3}
\end{equation}%
where $p_{0}^{\ast \ast }(x)$ and $g_{0,2}^{\ast \ast }(x,t)$ are functions
of our choice. The initial boundary value problem (\ref{8.1})-(\ref{8.3}) is
a conventional one and we solve it by the Finite Difference Method. Next, we
choose the coefficient $k\left( x\right) ,$ which we want to reconstruct.

Given the pair $\left( p^{\ast \ast }(x,t),k\left( x\right) \right) ,$ we
now need to ensure that the function ${v^{\ast \ast }(x,t)}$ satisfies the
first equation (\ref{2.5}). To do this, we define the function $s(x,t).$
Assume that $p^{\ast \ast }(x,t)\neq 0$ in $\overline{Q}_{T}.$ This can
always be achieved by a proper choice of functions $p_{0}^{\ast \ast }(x)$
and $g_{0,1}^{\ast \ast }(x,t)$ in (\ref{8.2}), (\ref{8.3}) and the use of
the maximum principle for parabolic equations. Then we define the function $%
s(x,t)$ as: 
\begin{eqnarray}
\hspace{-1 cm} \left. s(x,t)=\left[ p^{\ast \ast }(x,t)\right] ^{-1}\left[
v_{t}^{\ast \ast }+\Delta v^{\ast \ast }{-r(\nabla v^{\ast \ast })^{2}/2-}%
k\left( x\right) \int\limits_{\Omega }K\left( x,y\right) p^{\ast \ast
}\left( y,t\right) dy \right] ,\right.  \nonumber
\end{eqnarray}%
see (\ref{2.4}) for the integral operator. Thus, this is our numerical
method for the generation of the solution $\left( v^{\ast \ast
}(x,t),p^{\ast \ast }(x,t)\right) $ of MFGS (\ref{2.5}). Next, we use the
values of these two functions and their $\partial _{x_{1}}-$derivatives to
generate the target data (\ref{2.400}) for our CIP.

\subsection{Numerical arrangements}

\label{sec:8.2}

We have conducted numerical studies in the 2D case. In our numerical
testing, we took in (\ref{2.1}), (\ref{2.5}), (\ref{3.9}), (\ref{3.10}) and (%
\ref{3.11}): 
\[
a=1,\quad b=2,\quad A_{2}=0.5,\quad T=1,\quad r\left( x\right) \equiv
1,\quad \alpha =1/5,\quad \beta =0.001. 
\]%
Also, we have taken in (\ref{3.11}) 
\begin{equation}
\beta \left\Vert \left( u,m\right) \right\Vert _{\widetilde{H}}^{2}=\beta
\left( \left\Vert u\right\Vert _{H^{2}\left( Q_{T}\right) }^{2}+\left\Vert
m\right\Vert _{H^{2}\left( Q_{T}\right) }^{2}\right) .  \label{8.30}
\end{equation}

As to the kernel $K\left( x,y\right) $ in the global interaction term, we
have taken it as in (\ref{2.2}), 
\begin{equation}
K_{1}\left( x,y\right) =k\left( x\right) \delta \left( x_{1}-y_{1}\right) 
\overline{K}_{1}\left( x,y_{2}\right) ,  \label{8.4}
\end{equation}%
where $\overline{K}_{1}\left( x,y_{2}\right) $ is the Gaussian \cite[%
subsection 4.2]{LiuOsher}, 
\[
\overline{K}_{1}\left( x,y_{2}\right) =\exp \left( -\frac{\left(
x_{2}-y_{2}\right) ^{2}}{\sigma ^{2}}\right) ,\quad \sigma =0.2. 
\]%
We took the target coefficient $k(x)$ in (\ref{8.4}), to be reconstructed,
as: 
\begin{equation}
k(x)=\left\{ 
\begin{array}{cc}
c_{a}=const.>0, & \mbox{inside the tested inclusion,} \\ 
1, & \mbox{outside the tested inclusion.}%
\end{array}%
\right.  \label{8.5}
\end{equation}%
Then we set: 
\begin{equation}
\mbox{correct inclusion/background contrast}=\frac{c_{a}}{1}.  \label{8.6}
\end{equation}%
In the numerical tests below, we take $c_{a}=2,4,8$, and the inclusions with
the shapes of the letters `$A$', `$\Omega $' and `$SZ$'. To generate the
observation data in (\ref{2.400}), we set in (\ref{8.1})-(\ref{8.3}): 
\[
\hspace{-2cm}v^{\ast \ast }(x,t)=0.1\cos (\pi x_{1})\sin (\pi
x_{2})(t^{2}+1),\quad p_{0}^{\ast \ast }(x)=x_{1}x_{2},\quad g_{0,2}^{\ast
\ast }\left( x,t\right) =(t+1)x_{1}x_{2}. 
\]%
Then we proceeded with data generation as in subsection 8.1.

\textbf{Remark 8.1.} \emph{We choose letter-like shapes of inclusions
because these are non-convex shapes with voids. Shapes with such properties
are traditionally hard to image when solving inverse problems. Therefore,
our results point towards the robustness of our numerical technique.
Furthermore, we image both: shapes of inclusions and values of the unknown
coefficient }$k\left( x\right) $\emph{, }$x\in \Omega $\emph{, see (\ref{8.5}%
), (\ref{8.6}).}

To solve the forward problem (\ref{8.1})-(\ref{8.3}) for data generation, we
have used the spatial mesh sizes $1/80\times 1/80$ and the temporal mesh
step size $1/320$. In the computations of the Minimization Problem, the
spatial mesh sizes were $1/20\times 1/20$ and the temporal mesh step size
was $1/10$. We have solved problem (\ref{8.1})-(\ref{8.3}) by the classic
implicit scheme. To solve the Minimization Problem, we have written
operators $L_{1},L_{2}$ and the norm $\left\Vert \cdot \right\Vert _{%
\widetilde{H}}$ in (\ref{8.30}) in discrete forms of finite differences and
then minimized the functional $J_{\lambda ,\beta }$ with respect to the
values of functions $u$ and $m$ at those grid points. As soon as its
minimizer $\left( u_{\min },m_{\min }\right) $ is found, the computed target
coefficient $k_{\mbox{comp}}\left( x\right) $ is found via an obvious analog
of (\ref{3.5}), (\ref{3.6}).

To guarantee that the solution of the problem of the minimization of the
functional $J_{\lambda ,\beta }\left( u,m\right) $ in (\ref{3.11}) satisfies
the boundary conditions (\ref{3.63}), we adopt the Matlab's built-in
optimization toolbox \textbf{fmincon} to minimize the discretized form of
the functional $J_{\lambda ,\beta }\left( u,m\right) $. The iterations of 
\textbf{fmincon} stop when the condition 
\[
|\nabla J_{\lambda ,\beta }\left( u,m\right) |<10^{-2} 
\]
is met.

To exhibit the process for dealing with the Neumann boundary conditions in
the second line of (\ref{3.63}) in the iterations of \textbf{fmincon}, we
denote the discrete points along $x_{1}$-direction as 
\begin{equation}
\hspace{-1cm}x_{1,i}=a+ih_{x_{1}},\quad i=0,1,\cdots ,N_{x_{1}},,\quad
N_{x_{1}}=(b-a)/h_{x_{1}},\quad h_{x_{1}}=1/20.  \label{8.801}
\end{equation}%
Then keepng in mind those Neumann boundary conditions, the discrete
functions $u,m$ in the iterations of \textbf{fmincon} should satisfy 
\begin{equation}
\hspace{-1.8cm}\left. 
\begin{array}{c}
-4u(x_{1,N_{x_{1}-1}},x_{2},t)+u(x_{1,N_{x_{1}-2}},x_{2},t)=2h_{x_{1}}%
\partial _{t}g_{0,1}(b,x_{2},t)-3\partial _{t}g_{1,1}(b,x_{2},t), \\ 
-4m(x_{1,N_{x_{1}-1}},x_{2},t)+m(x_{1,N_{x_{1}-2}},x_{2},t)=2h_{x_{1}}%
\partial _{t}g_{0,2}(b,x_{2},t)-3\partial _{t}g_{1,2}(b,x_{2},t).%
\end{array}%
\right.   \label{8.802}
\end{equation}%
We note that formula \eref{8.802} also contains the Dirichlet boundary
conditions in the first line of (\ref{3.63}) as 
\begin{equation}
\left. 
\begin{array}{c}
u(x_{1},x_{2},t)=\partial _{t}g_{0,1}(x_{1},x_{2},t),\mbox{ }%
m(x_{1},x_{2},t)=\partial _{t}g_{0,2}(x_{1},x_{2},t),\mbox{ }x\in S_{T}.%
\end{array}%
\right.   \label{8.803}
\end{equation}

The starting point $\left( u^{\left( 0\right) },m^{\left( 0\right) }\right)
\left( x_{1},x_{2},t\right) $ of iterations of \textbf{fmincon} was chosen
as: 
\begin{equation}
\hspace{-1.5cm}%
\eqalign{ u^{(0)}(x_{1}, x_{2}, t)= &\frac{1}{2} \left(
\frac{b-x_{1}}{b-a} \partial_{t}g_{0,1}\left(x_{1}, x_{2}, t\right) +
\frac{x_{1}-a}{b-a} \partial_{t}g_{0,1}\left(x_{1}, x_{2}, t\right) \right)
\\ &+ \frac{1}{2} \left( \frac{A_{2}-x_{2}}{2A_{2}}
\partial_{t}g_{0,1}\left(x_{1}, x_{2}, t\right) + \frac{x_{2}+A_{2}}{2A_{2}}
\partial_{t}g_{0,1}\left(x_{1}, x_{2}, t\right) \right), \\ m^{(0)}(x_{1},
x_{2}, t)= &\frac{1}{2} \left( \frac{b-x_{2}}{b-a}
\partial_{t}g_{0,2}\left(x_{1}, x_{2}, t\right) + \frac{x_{2}-a}{b-a}
\partial_{t}g_{0,2}\left(x_{1}, x_{2}, t\right) \right) \\ &+ \frac{1}{2}
\left( \frac{A_{2}-x_{2}}{2A_{2}} \partial_{t}g_{0,2}\left(x_{1}, x_{2},
t\right) + \frac{x_{2}+A_{2}}{2A_{2}} \partial_{t}g_{0,2}\left(x_{1}, x_{2},
t\right) \right). }  \label{8.804}
\end{equation}%
In other words, we use linear interpolation inside the domain $Q_{T}$ of the
Dirichlet boundary conditions (\ref{3.63}). Although it follows from (%
\ref{8.804}) that the starting point $\left( u^{\left( 0\right) },m^{\left(
0\right) }\right) \left( x_{1},x_{2},t\right) $ satisfies only Dirichlet
boundary conditions in the first line of (\ref{3.63}) and does not satisfy
Neumann boundary conditions in the second line of (\ref{3.63}), still (\ref%
{8.802}) and (\ref{8.803}) imply that boundary conditions in both lines of (%
\ref{3.63}) are satisfied on all other  iterations of \textbf{fmincon}.

We introduce the random noise in the observation data in (\ref{2.400}) as
follows: 
\begin{equation}
\hspace{-1.5cm}\eqalign{ & \hspace{1.6cm}v_{0,\zeta }(x)=v_{0}(x)\left(
1+\delta \zeta _{v,x}\right) ,\quad p_{0,\zeta }(x)=p_{0}(x)\left( 1+\delta
\zeta _{p,x}\right) , \\ & g_{0,1,\zeta }(x,t)=g_{0,1}(x,t)\left( 1+\delta
\zeta _{0,1,x,t}\right) ,\quad g_{0,2,\zeta }(x,t)=g_{0,2}(x,t)\left(
1+\delta \zeta _{0,2,x,t}\right) , \\ & g_{1,1,\zeta
}(x,t)=g_{1,1}(x,t)\left( 1+\delta \zeta _{1,1,x,t}\right) ,\quad
g_{1,2,\zeta }(x,t)=g_{1,2}(x,t)\left( 1+\delta \zeta _{1,2,x,t}\right) , }
\label{8.7}
\end{equation}%
where $\zeta _{v,x},\zeta _{p,x}$ are the uniformly distributed random
variables in the interval $[0,1]$ depending on the point $x\in \Omega $, $%
\zeta _{0,1,x,t},\zeta _{0,2,x,t}$ are the uniformly distributed random
variables in the interval $[0,1]$ depending on the point $(x,t)\in S_{T}$,
and $\zeta _{1,1,x,t},\zeta _{1,2,x,t}$ are the uniformly distributed random
variables in the interval $[0,1]$ depending on the point $(x,t)\in \Gamma
_{T}$. In (\ref{8.7}) $\delta =0.03,0.05$, which correspond to the $3\%$ and 
$5\%$ noise levels respectively. The reconstruction from the noisy data is
denoted as $k_{\zeta }(x)$ given by the following analog of (\ref{3.5}), (%
\ref{3.6}): 
\begin{equation}
\left. 
\begin{array}{c}
k_{\zeta }(x)=\left( u\left( x,t\right) -\int\limits_{T/2}^{t}u_{t}\left(
x,\tau \right) d\tau \right) f_{\zeta }(x)+F_{\zeta }(x), \\ 
f_{\zeta }(x)=\left[ \int\limits_{\Omega _{1}}\overline{K}_{1}\left( x,%
\overline{y}\right) p_{0,\zeta }\left( x_{1},\overline{y}\right) d\overline{y%
}\right] ^{-1},\mbox{ } \\ 
F_{\zeta }\left( x\right) =\left( \Delta v_{0,\zeta }-r\left( \nabla
v_{0,\zeta }\right) ^{2}/2-s\left( x,T/2\right) p_{0,\zeta }\right) \left(
x\right) f_{\zeta }\left( x\right) ,%
\end{array}%
\right.  \label{8.8}
\end{equation}%
where the subscript $\zeta $ means that these functions correspond to the
noisy data. Since we deal with first and second derivatives of noisy
functions $v_{0,\zeta }(x)$, $p_{0,\zeta }(x),g_{0,1,\zeta }(x,t)$, $%
g_{0,2,\zeta }(x,t)$, $g_{1,1,\zeta }(x,t)$, $g_{1,2,\zeta }(x,t)$, we have
to design a numerical method to differentiate the noisy data. We use the
natural cubic splines to approximate the noisy observation data. Then we can
use the derivatives of those splines to approximate the derivatives of
corresponding noisy observation data. We generate the cubic splines $%
s_{v},s_{p}$ of functions $v_{0,\zeta }(x),p_{0,\zeta }(x)$ with the spatial
mesh grid size of $1/20\times 1/20$, and then calculate the first and second
derivatives of $s_{v}$ to approximate the first and second derivatives
respect to $x,y$ of $v_{0,\zeta }(x),p_{0,\zeta }(x)$. For the boundary data 
$g_{i,j,\zeta }(x,t),i=0,1,j=1,2$, we generate the corresponding cubic
splines $S_{i,j},i=0,1,j=1,2$ in the temporal space with the temporal mesh
grid size of $1/10$, and then calculate the derivatives of $S_{i,j}$ to
approximate the first derivatives with respect to $t$ of functions $%
g_{i,j,\zeta }(x,t)$.

\subsection{Numerical tests}

\label{sec:8.3}

\begin{figure}[tbph]
\centering
\includegraphics[width = 6in]{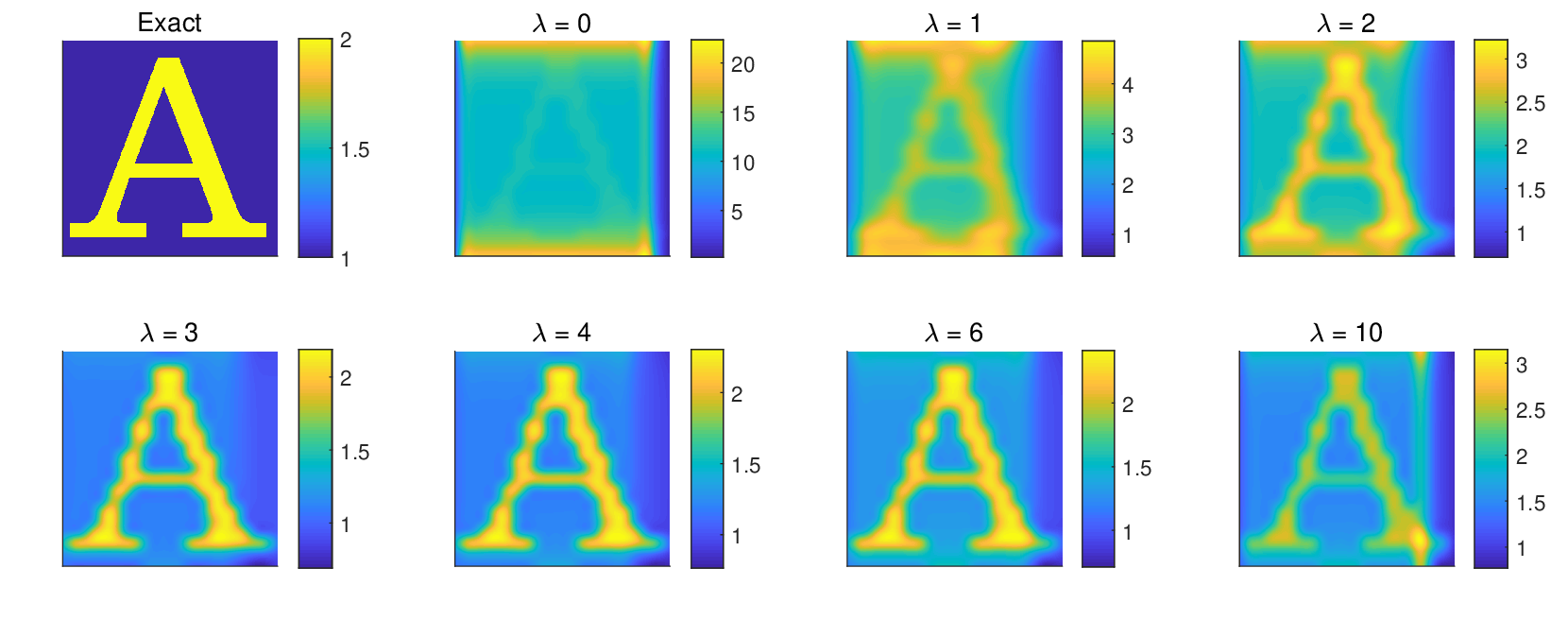}
\caption{Test 1. The reconstructed coefficient $k(x)$, where the function $%
k(x)$ is given in \eref{8.5} with $c_{a}=2$ inside of the letter `A'. We
test different values of $\protect\lambda $. The reconstructions are
unsatisfactory when $\protect\lambda $ is too small, $\protect\lambda =0,1,2$%
. Next, the reconstruction quality is good when $\protect\lambda =3,4$. On
the other hand, the reconstruction quality deteriorates at $\protect\lambda %
=10$. Thus, we choose $\protect\lambda =3$ as the optimal value. }
\label{plot_A_diff_lambda}
\end{figure}

\textbf{Test 1.} We test the case when the inclusion in (\ref{8.5}) has the
shape of the letter `$A$' with $c_{a}=2$. We use this test as a reference
case to figure out an optimal value of the parameter $\lambda $. The result
is displayed in Figure \ref{plot_A_diff_lambda}. We observe that the images
have a low quality for $\lambda =0,1,2$. Then the quality is improved with $%
\lambda =3,4$, and the reconstruction quality deteriorates at $\lambda =10$.
Hence, we choose $\lambda =3$ as the optimal value.

\textbf{Remark 8.2.}\emph{\ The optimal value }$\lambda =3,$\emph{\ once
chosen, is used in other tests 2-4.}

\begin{figure}[tbph]
\centering
\includegraphics[width = 4in]{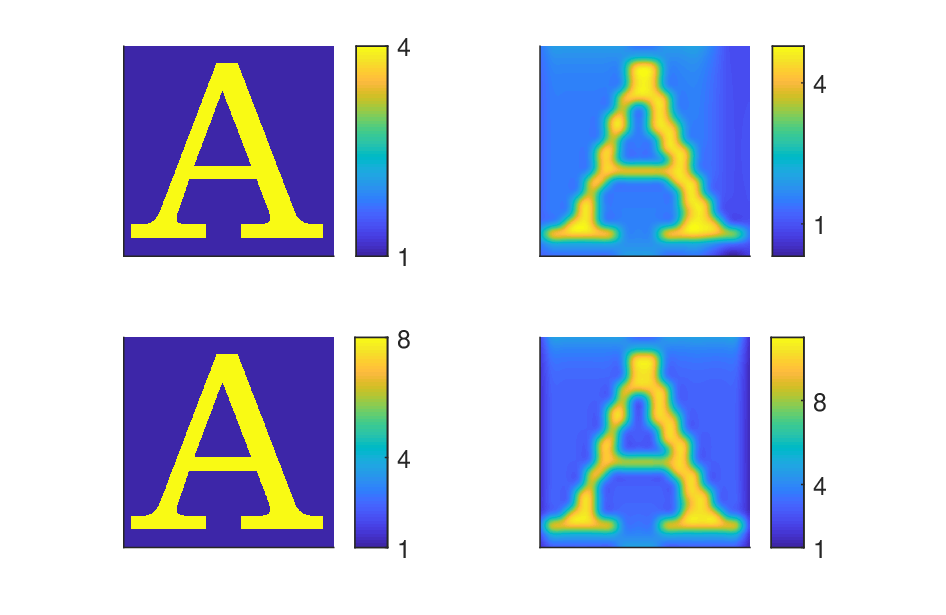}
\caption{Test 2. Exact (left) and reconstructed (right) coefficient $k(x)$
with $c_{a}=4$ (first row) and $c_{a}=8$ (second row) inside of the letter
`A' as in (\protect\ref{8.5}). The inclusion/background contrasts in (%
\protect\ref{8.6}) are respectively $4:1$ and $8:1$. The reconstructions are
accurate}
\label{plot_A_4_8}
\end{figure}

\textbf{Test 2.} We test the case when the inclusion in (\ref{8.5}) has the
shape of the letter `$A$' for different values of the parameter $c_{a}=4,8$
inside of the letter `$A$'. Hence, by (\ref{8.6}) the inclusion/background
contrasts now are respectively $4:1$ and $8:1$. Computational results are
displayed on Figure \ref{plot_A_4_8}. One can observe that these images are
accurate ones. In particular, the computed inclusion/background contrasts
are accurate. 

\begin{figure}[tbph]
\centering
\includegraphics[width = 4in]{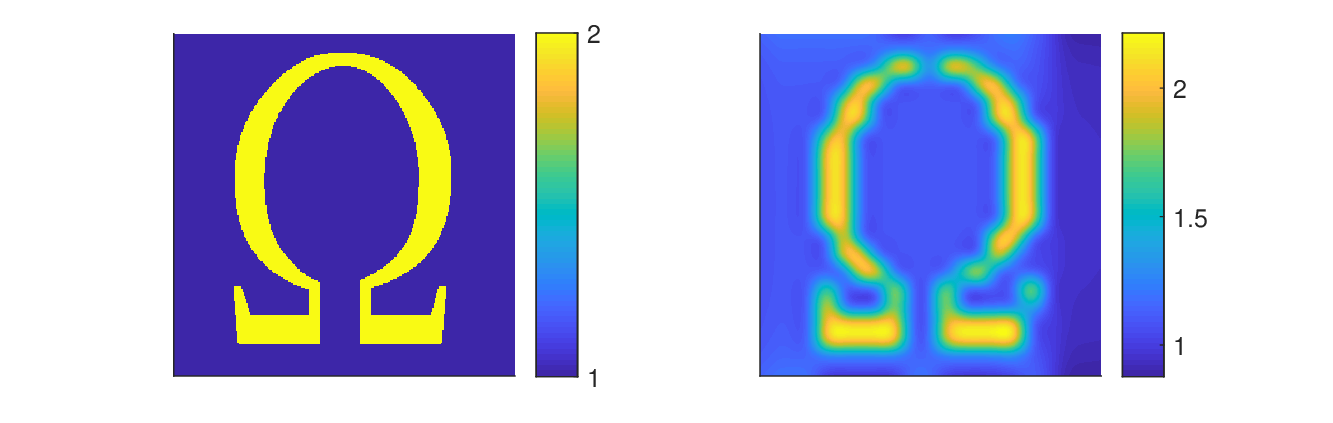}
\caption{Test 3. Exact (left) and reconstructed (right) coefficient $k(x)$,
where the function $k(x)$ is given in \eref{8.5} with $c_{a}=2$ inside of
the letter `$\Omega $'. The reconstruction is accurate.}
\label{plot_Omega}
\end{figure}

\textbf{Test 3.} We test the case when the coefficient $k(x)$ in \eref{8.5}
has the shape of the letter `$\Omega $' with $c_{a}=2$ inside of it. Results
are presented on Figure \ref{plot_Omega}. We again observe an accurate
reconstruction. 

\begin{figure}[htbp]
\centering
\includegraphics[width = 4in]{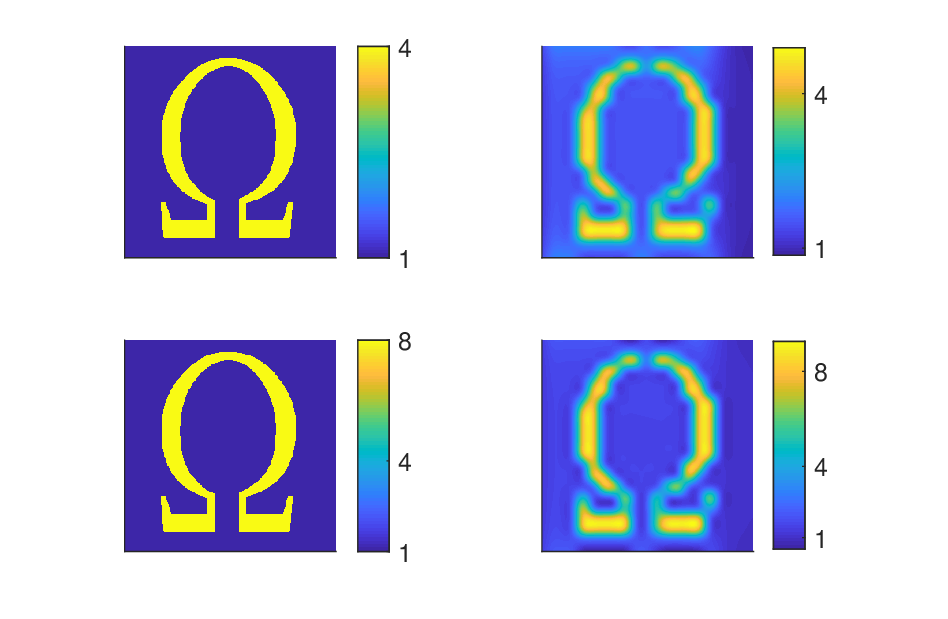}
\caption{Test 4. Exact (left) and reconstructed (right) coefficient $k(x)$
with $c_{a}=4$ (first row) and $c_{a}=8$ (second row) inside of the letter `$%
\Omega$' as in (\protect\ref{8.5}). The inclusion/background contrasts in (%
\protect\ref{8.6}) are respectively $4:1$ and $8:1$. The reconstructions are
accurate.}
\label{plot_Omega_4_8}
\end{figure}

\textbf{Test 4.} We now test different values of the parameter $c_{a}=4,8$
inside of the letter `$\Omega $'. The computational results are depicted on
Figure \ref{plot_Omega_4_8}. The quality of these images and computed
inclusion/background contrasts in (\ref{8.6}) are as good as the results of
the test with the letter `$A$' in Test 2.

\textbf{Test 5.} We test the reconstruction for the case when the inclusion
in (\ref{8.5}) has the shape of two letters `SZ' with $c_{a}=2$ in each of
them. SZ are two letters in the name of the city (Shenzhen) were the second
and the third authors reside. The results are exhibited on Figure \ref%
{plot_SZ}. 
\begin{figure}[tbph]
\centering
\includegraphics[width = 4in]{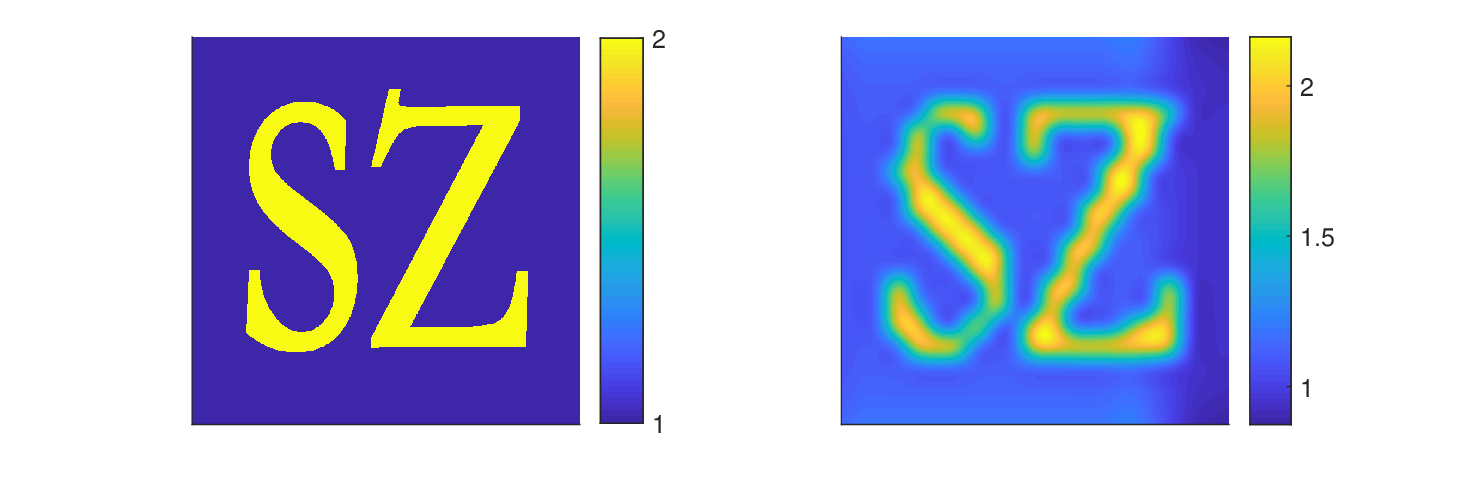}
\caption{Test 5. Exact (left) and reconstructed (right) coefficient $k(x)$,
where the function $k(x)$ is given in \eref{8.5} with $c_{a}=2$ inside of
two letters `SZ'. The reconstruction is worse than the one for the case of
the single letter `$\Omega $' in Figure \protect\ref{plot_Omega}.
Nevertheless, the reconstruction are still accurate in both letters. The
computed inclusion/background contrasts in \eref{8.6} are accurately
reconstructed in both letters.}
\label{plot_SZ}
\end{figure}

\begin{figure}[tbph]
\centering
\includegraphics[width = 4in]{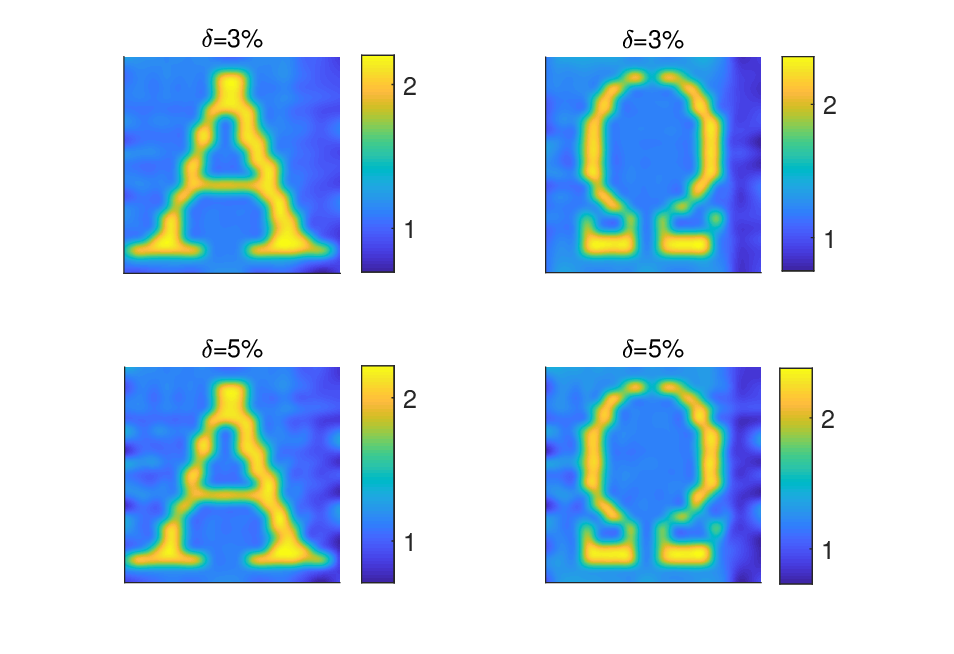}
\caption{Test 6. Reconstructed coefficient $k(x)$ with the shape of letters
`A' (left) and `$\Omega $' (right) with $c_{a}=2$ from noisy data \eref{8.7}
with $\protect\delta =0.03$ (top) and $\protect\delta =0.05$ (bottom), i.e.
with 3\% and 5\% noise level. In all these four cases, both reconstructions
and inclusion/background contrasts in \eref{8.6} are accurate.}
\label{plot_Noise}
\end{figure}

\textbf{Test 6.} We consider the case of the random noisy data in (\ref{8.7}%
) with $\delta =0.03$ and $\delta =0.05$, i.e. with 3\% and 5\% noise level.
We test the reconstruction for the cases when the inclusion in (\ref{8.5})
has the shape of either the letter `A' or the letter `$\Omega $' with $%
c_{a}=2$. The results are displayed on Figure \ref{plot_Noise}. One can
observe accurate reconstructions in all four cases. In particular, the
inclusion/background contrasts in (\ref{8.6}) are reconstructed accurately.

\section{Summary}

\label{sec:9}

For the first time, we have developed here a numerical method for a CIP for
MFGS (\ref{2.5}) with the rigorously guaranteed global convergence property,
see Definition 1.2 for the term \textquotedblleft global convergence". In
other words, its convergence to the true solution is guaranteed if starting
at an arbitrary point of a bounded set, whose diameter is not required to be
small. We have proven the global convergence of our technique and explicitly
provided its convergence rate.

To computationally demonstrate robustness of our technique, we took quite
complicated shapes of tested inclusions. Our numerical results for both
noiseless and noisy data demonstrate that we can accurately reconstruct both
shapes of inclusions and values of the unknown coefficient both inside and
outside of them.

\section*{References}


\end{document}